%% file: MAKPPAS-kmc.tex
\documentclass{siamltex}

\usepackage{amsfonts,amssymb}        
\usepackage{amsmath}
\usepackage{epsfig}

\usepackage{times}

\def\MAKCOMM#1{}
\def\ALEXCOMM#1{}

\input{ppmacro1.tex}

\title{%
Error analysis of coarse-grained kinetic Monte Carlo method
}

\author{%
Markos A. Katsoulakis\thanks{%
Department of Mathematics and Statistics, University of Massachusetts,
Amherst, MA 01003--9305, USA,
{\tt markos@math.umass.edu}}
\and
Petr Plech\'a\v{c}\thanks{%
Mathematics Institute, The University of Warwick,
Coventry, CV4 7AL, United Kingdom, 
{\tt plechac@maths.warwick.ac.uk}}
\and
Alexandros Sopasakis\thanks{%
Department of Mathematics and Statistics, University of Massachusetts,
Amherst, MA 01003--9305, USA,
{\tt sopas@math.umass.edu}}
}
\begin{document}

\maketitle

\begin{abstract}
The coarse-grained Monte Carlo (CGMC) algorithm was originally proposed
in the series of works \cite{kmv1,kmv2, kv}. 
In this paper we further investigate the approximation properties
of the coarse-graining procedure and provide both analytical and numerical 
evidence that the hierarchy of the coarse models is built in a systematic
way that allows for  error control 
in both transient and long-time simulations. We demonstrate that the numerical 
accuracy of the CGMC algorithm as an approximation of stochastic lattice spin 
flip dynamics is of  order two  in terms of the coarse-graining ratio and that
the natural small parameter is the coarse-graining ratio  
over the range of particle/particle interactions. The error estimate is shown 
to hold   in the weak convergence sense.
We employ the 
derived analytical results to guide CGMC algorithms
and we demonstrate a CPU speed-up in demanding computational  regimes that 
involve nucleation, phase transitions  and   metastability.
\end{abstract}


%
%

%
%
\section{Introduction} \label{sec:level1}

%
Microscopic  computational models for complex systems  such as Mo\-lecular Dynamics (MD) and Monte Carlo 
(MC) algorithms are typically  formulated in terms of simple rules describing interactions between 
individual particles or spin variables. 
The large number of variables 
and even larger number of interactions between them present the principal limitation for 
efficient simulations. 
     Another restricting factor is illustrated by essentially sequential nature
     of approximating the time evolution in  particle systems that yields a substantial
        slowdown in the resolution of dynamics, especially in metastable regimes.


In \cite{kmv1, kmv2, kv} the authors started developing 
systematic mathematical strategies for the coarse-graining of microscopic models, focusing on the paradigm 
of stochastic  lattice dynamics and the corresponding MC simulators. In principle, coarse-grained models 
are expected to have fewer observables than the original microscopic system making them computationally 
more efficient than the direct numerical simulations. In these papers  a   hierarchy of coarse-grained 
stochastic models -- referred to as coarse-grained MC (CGMC) -- was derived from the microscopic rules
through  a stochastic closure argument.    
The CGMC hierarchy  is reminiscent of Multi-Resolution Analysis 
approaches to the discretization of operators \cite{bcr}, spanning length/time scales from the microscopic 
to the mesoscopic. The resulting  {\it stochastic coarse-grained processes}  involve Markovian birth-death 
and generalised exclusion processes and their combinations, and as  demonstrated in \cite{kmv1,kmv2,kv},  
they share the same ergodic properties with their microscopic counterparts.  
The full hierarchy of the coarse-grained stochastic dynamics satisfies detailed balance relations 
and as a result not only it yields self-consistent random fluctuation mechanism, but also consistent 
with the underlying microscopic fluctuations and the unresolved degrees of freedom.  
From the computational complexity perspective, a comparison of CGMC with conventional MC methods for 
the same real time shows, \cite{kmv1}, that the CPU time can decrease approximately as $O(1/q^2)$ or faster,
where $q$ is the level of coarse-graining, 
as demonstrated for spin-flip lattice dynamics. 
Thus, while for 
macroscopic size systems in the millimeter length scale or larger, microscopic MC simulations are 
impractical on a single processor, the computational savings of CGMC make it a suitable tool 
capable of capturing large scale features, while retaining microscopic information on intermolecular 
forces and particle fluctuations. In the case of diffusion (spin exchange) dynamics 
we also observe 
an additional coarse-graining in time by a factor $q^2$, improving the hydrodynamic slowdown effect 
in the conservative MC, \cite{kv}.  


In the recent paper \cite{kt}
the authors  rigorously analysed  CGMC models as 
approximations of conventional MC in  {\it non-equilibrium},  by estimating the {\it information loss} 
between microscopic and coarse-grained adsorp\-tion/desorpt\-ion lattice dynamics. In analogy to  
the numerical analysis 
for PDEs,  an error 
analysis was carried out  between   the {\it exact microscopic   
process} $\PROCMIC$ and the  {\it approximating coarse-grained process}  $\PROCMAC$. 
The key  step in this direction was to use, as  a quantitative measure for    the   
{loss of information} in the  coarse-graining from finer to coarser scales, the 
information-theoretic concept of   the {\it relative entropy} between probability measures, \cite{ct}. 
Such relative entropy estimates give a first mathematical reasoning for the parameter regimes, i.e., 
the  degree of coarse-graining versus the interaction range, for which CGMC is expected to give errors 
within a given tolerance.  
Using the rigorous results in \cite{kt} as a starting point, in this paper  we focus on carrying 
out a detailed numerical analysis of the error propagation for spin flip lattice dynamics. 
{Due to the numerical intractability of the  relative entropy 
for a large particle system, we employ, in the numerical error calculations,  suitable computable 
upper and lower bounds, as well as {\it  targeted } coarse observables.
The latter point of view  
necessitates in the use of a weak convergence framework for the study of the error
between CGMC and direct numerical simulations of the stochastic lattice dynamics.
We demonstrate that the numerical accuracy of the CGMC algorithm  is of  order two  in terms of 
the ratio of the coarse-graining over the range of particle/particle interactions. 
We also refer to recent work in \cite{KS} on weak error estimates
between microscopic MC algorithms and therein derived  SDE  approximations.
Further details about a priori estimates for weak convergence of
approximations to SDEs can be found in \cite{Talay-Tubaro1,Talay-Tubaro2} 
and \cite{Kloeden}. Related a posteriori estimates are discussed in
\cite{Szepessy1}.
We further employ the 
derived analytical results to guide CGMC algorithms
and we demonstrate a CPU speed-up in demanding computational  regimes that involve nucleation, 
phase transitions  and   metastability.
We demonstrate computationally that CGMC probes efficiently the energy landscape, 
yielding {\it  spatial path-wise } agreement with the underlying microscopic lattice dynamics, 
at least for fairly long but still finite interactions. }

The CGMC algorithms discussed here are related to a number of methods involving 
coarse-graining at various levels, for instance fast summation techniques, computational
re\-nor\-malization and simulation  and multi-scale computational methods for stochastic systems. 
One of the sources of the computational  complexity of molecular simulations arises in the 
calculation of particle/particle interactions, especially in the case where long range 
forces are  relevant. The evaluation cost of  such pairwise interactions can be significantly 
reduced by applying well-controlled approximation schemes and/or a hierarchical decomposition 
of the computation. 
Such ideas have been successfully applied in the development of Ewald 
summation techniques,  multigrid  (MG) ,  fast multipole methods (FMP) or tree-code algorithms.
Typically, once the interaction terms are computed with one of these fast summation methods, 
they are entered in the microscopic algorithm where  a simulation with a large number of 
individually tracked particles has still to be carried out. The point of view adopted by 
CGMC is related to these methods in the sense that the interaction potential or 
operator is approximated  in terms of a  truncated multi-resolution decomposition  
within a given tolerance. The CGMC is subsequently defined at the coarse level specified 
by the truncation of the decomposition. However, a notable difference is that CGMC models 
track much fewer  coarse observables instead of simulating  every  individual particle.  
The equilibrium set-up of CGMC  is essentially given by  the renormalised Hamiltonian after 
a single  iteration in the renormalisation group flow. It is not surprising that such 
an approach, when applied to near critical temperature simulations,
has many limitations.  For example, in  the nearest-neighbour Ising-type models 
this fact is manifested  in the aforementioned  
error estimates and the  comparative simulations in \cite{kmv1}.  On the other hand the 
focus of CGMC is dynamic simulations usually  coupled to a macroscopic system (see for 
instance the hybrid systems in \cite{vlachos90, kms}), where criticality may not be as 
important due to the presence of a time-varying external field. Nevertheless, further  
corrections to the CGMC dynamics from the renormalisation group
flow given by RGMC and multigrid MC methods
\cite{B92,BrandtRon:01,Sokal} can improve the order of convergence  of the CGMC.
We refer to \cite{KPRT} for higher order accurate CGMC methods based on cluster expansions, 
where the coarse-graining procedure described here is the model  around which a cluster 
expansion is carried out with controlled errors. 
In that sense the CGMC method 
is of order two  accurate as explained in Section 4.

In recent years there has been a growing interest in developing and analysing coarse-graining methods  
for the purpose of modelling and simulation across scales. Such systems arise   in a broad 
spectrum of scientific disciplines ranging from materials science to macromolecular dynamics, 
to epidemiology  and to atmosphere/ocean science. 
Various coarse-graining approaches may yield explicitly 
derived stochastic coarse models using different coarse approximations, e.g.,
\cite{golden, steph, steph1,mtv,schuette1}, or 
can be statistics-based \cite{m-p} or may rely on on-fly simulations, e.g., the equation-free method 
\cite{kevrek}, the heterogeneous multi-scale method  \cite{ee}, or multi-scale FE methods \cite{hou}. 
A systematic approach to upscaling of stochastic systems has been also developed from the
multi-level perspective
in \cite{BRA,Bai,BI2}, where the authors proposed algorithms
for efficient multi-scale simulations using Monte Carlo methods.
Other coarse-graining
techniques in the polymer science literature include the bond fluctuation model and its 
variants \cite{bi}.
Such coarse-graining methodologies often rely on parametrisation,
hence at different conditions (e.g., temperature, density, composition) coarse potentials 
need to be re-parametrised \cite{m-p}. 
%

%
%
\section{Microscopic lattice models}

The presented analysis applies to the class of Ising-type lattice systems.
For the sake of simplicity we assume that the computational domain is
defined as the discrete periodic lattice $\LATT = \frac{1}{n} \Z^d \cap \T$
which represents discretion of the $d$-dimensional torus $\T=[0,1)^d$ 
and  $d$ denotes the spatial dimension. We restrict presentation of the
results to $d=1$, nevertheless higher dimensional cases are obtained without
significant changes. 
However, the algorithms can also be implemented
on bounded domains with usual boundary conditions. 
The number of lattice sites $N=n^d$ is fixed. 
The microscopic degrees of freedom or
the microscopic order parameter is given by the spin-like 
variable $\sigma(x)$ defined at each
site $x\in\LATT$. In this paper we discuss only the case of discrete 
spin variables, i.e., $\sigma(x)\in \SPINSP$ with $\SPINSP=\{-1,1\}$,
$\SPINSP=\{0,1\}$ (Ising model) or $\SPINSP=\{0,1,\dots s\}$ (Potts models).
The case of the spin variable belonging to a compact Riemannien
manifold, e.g., $\SPINSP= \Sph^2$ (Heisenberg model), 
$\SPINSP=\mathrm{SU}(2)$ (matrix model), will be studied elsewhere.
We denote by $\CONFMICRO=\{\sigma(x)\SEP x\in\LATT\}$ a configuration of spins on
the lattice, i.e., an element of the configuration space $\SIGMA = \SPINSP^{\LATT}$.
The interactions between spins at a given configuration $\CONFMICRO$
are defined by the microscopic Hamiltonian
\begin{equation}\label{hamiltonian}
H(\CONFMICRO) = -\frac{1}{2}\sum_{x\in\LATT}\sum_{y\neq x} 
                  J(x-y) \sigma(x)\sigma(y) + \sum_{x\in\LATT} h(x)\sigma(x)\COMMA
\end{equation}
where $h(x)$ denotes the external field at the site $x$. The two-body inter-particle
potential $J$ accounts for interactions between individual spins. We consider
the class of potentials with the following properties
\begin{eqnarray}
&& J(x-y) = \frac{1}{L^d}\JNOT\left(\frac{n}{L}|x-y|\right)\COMMA
   \SPACE x,y\in\LATT\COMMA \label{defJV1}\\
&& \JNOT:\R\to\R\COMMA\SPACE \JNOT(r)=\JNOT(-r)\COMMA\SPACE
   \JNOT(r) = 0 \COMMA\SPACE\mbox{ if $|r|\geq 1$.} \label{defJV2} 
\end{eqnarray}
We impose additional assumptions on $\JNOT$ which allow us to derive
explicit error estimates:
\begin{eqnarray}
&& \mbox{$\JNOT$ is smooth on $\R\setminus\{0\}$,}\label{AssumpV1} \\
&& \int_{\R}|\JNOT(r)| \,dr < \infty\COMMA\mbox{ and}\; 
   \int_{\R} |\partial_r \JNOT(r)|\,dr < \infty\PERIOD \label{AssumpV2}
\end{eqnarray}
Note that the summability condition for $\JNOT$ guarantees that the potential
$J$ is also summable due to the scaling factor. Hence the Hamiltonian is well 
defined even for $N,L\to\infty$. The canonical equilibrium state is given
in terms of the Gibbs measure
\begin{eqnarray}\label{gibbs}
&&\mu_{N,\beta} (d\sigma) = \frac{1}{Z_{N,\beta}} 
             \EXP{-\beta H(\sigma)} \PRODM_N(d\sigma)\COMMA\SPACE\SPACE
  Z_{N,\beta} = \int_\SIGMA \EXP{-\beta H(\sigma)} \PRODM_N(d\sigma)\COMMA
\end{eqnarray}
where $\PRODM_N(d\sigma) = \PROD_{x\in\LATT} \rho(d\sigma(x))$ is the product
measure on $\SIGMA$ and the spins $\sigma(x)$ are independent identically
distributed (i.i.d.) random variables with the common distribution $\rho$.
For example, in  the Ising model the prior distribution on $\SPINSP=\{0,1\}$
would typically be $\rho(0) =\rho(1) = 1/2$. 

\smallskip

The microscopic dynamics is defined as a continuous-time jump  Markov process 
that defines a change of the spin  $\sigma(x)$ with the probability $c(x,\sigma;\xi)\DT$
over the time interval $[t,t+\DT]$. The function $c:\LATT\times \SIGMA\times\SPINSP\to\R$ is called
a rate of the process. The jump process $\PROCMIC$ is constructed in the following way:
suppose that at the time $t$ the configuration is $\PROCMICRO$, then the probability of
changing  the spin at the site $x\in\LATT$ spontaneously 
from $\PROCMICRO(x)$ to a new value $\xi\in\SPINSP$ over the time interval $[t,t+\DT]$ is 
$c(x,\sigma;\xi)\DT + O(\DT^2)$. We denote the resulting configuration by
$\sigma^{x,\xi}$. In the case of the Ising-type state space and spin-flip dynamics
we omit $\xi$ in this notation. 
The generator $\LOPER:L^\infty(\SIGMA)\to L^\infty(\SIGMA)$
of the Markov process acting on a bounded test function $\phi\in L^\infty(\SIGMA)$
defined on the space of configurations 
is given by
\begin{equation}\label{generator}
(\LOPER \phi)(\sigma) = \sum_{x\in\LATT} \int_{\SPINSP} c(x,\sigma;\xi) 
                  \left(\phi(\sigma^{x,\xi}) - \phi(\sigma)\right)\,d\xi\PERIOD
\end{equation}
The evolution of an observable (a test function) $\phi$  is given by 
\begin{equation}
\frac{d}{dt}\EXPEC{\phi(\PROCMICRO)} = \EXPEC{\LOPER \phi(\PROCMICRO)}\COMMA
\end{equation}
where the expectation operator $\EXPEC{.}$ is with respect to a measure conditioned
to the initial configuration $\sigma_{t=0} = \sigma_0$.
We require that the dynamics is of relaxation type such that the invariant measure
of this Markov process is the Gibbs measure \VIZ{gibbs}. The sufficient condition 
is known as {\it Detailed Balance} (DB)  and it imposes condition on the form
of the rate 
\begin{equation}\label{DB}
c(x,\sigma;\xi) \EXP{-\beta H(\sigma)} = 
        c(x,\sigma^{x,\xi};\sigma(x)) \EXP{-\beta H(\sigma^{x,\xi})}\PERIOD
\end{equation}
This condition has a simple interpretation: 
$c(x,\sigma;\xi)$ is the rate of converting $\sigma(x)$ to the value
$\xi$ while  $c(x,\sigma^{x,\xi};\sigma(x))$ is the rate of changing the spin
with the value $\xi$ at the site $x$ back to $\sigma(x)$.
The widely used class of Metropolis-type dynamics satisfies \VIZ{DB} and has the
rate given by
\begin{equation}
c(x,\sigma;\xi) = G(\beta \Delta_{x,\xi}H(\sigma))\COMMA\mbox{ where 
        $\Delta_{x,\xi}H(\sigma) = H(\sigma^{x,\xi}) - H(\sigma)$,}
\end{equation}
where $G$ is a continuous function satisfying: $G(r)=G(-r)e^{-r}$ for
all $r\in\R$. The most common choices in physics simulations are
$G(r)=\frac{1}{1+e^r}$ (Glauber dynamics), $G(r)=e^{-[r]_+}$, (Metropolis dynamics), with
$[r]_+ = r$ if $r\geq 0$ and $=0$ otherwise, or $G(r)=e^{-r/2}$. Such dynamics are often
used as samplers from the canonical equilibrium Gibbs measure. However, the
kinetic Monte Carlo method is also used for simulations of non-equilibrium processes.
The dynamics in such a case is known as {\it Arrhenius dynamics}, whose  rates
are usually derived from transition state theory or obtained from molecular dynamics 
simulations. 

To avoid unnecessary generality we restrict the description
to the Ising-type model with $\SPINSP=\{0,1\}$ used for modelling adsorption/desorption
processes. We also omit $\xi$ in the notation. The Arrhenius rate is defined
as follows
\begin{equation}\label{arrhenius_rate}
c(x,\sigma) = \left\{ \begin{array}{cl}
                      d_0    &  \mbox{if $\sigma(x)=0$,} \\
                     d_0\EXP{-\beta U(x,\sigma)} &\mbox{if $\sigma(x)=1$,}
                      \end{array}\right.
\end{equation}
where 
\begin{equation}\label{potentialU}
U(x,\sigma) = \sum_{y\in\LATT\,, y\neq x} J(x-y)\sigma(y) - h(x)\PERIOD
\end{equation}
Furthermore, the spin-flip rule is given by
$$
\sigma^{x}(y) = \left\{ \begin{array}{cl}
                      1 - \sigma(x)    &  \mbox{if $y=x$} \\
                      \sigma(y)        &  \mbox{if $y\neq x$.}
                        \end{array}\right.
$$
With the introduced notation the coarse-graining algorithm can be described
as an {\it approximation} of the microscopic dynamics, i.e., of the process 
$\PROCMIC$ by a coarse-grained process $\PROCMAC$ where the approximation
is done in a controlled way. We are interested not only in the approximation
of the invariant measure $\mu_{N,\beta} (d\sigma)$ (see \VIZ{gibbs}) but also in
the approximation of the measure on the path space.

%
%
\section{Approximation of the coarse-grained process}

The coarse-graining is defined in a geometric way by introducing the coarse-grained
observables as  block-spin variables. This approach follows the standard procedure of real-space
renormalisation, see for example \cite{Kadanoff}. We remark that although we introduce
block-spins our aim is not to approximate the renormalisation group  flow (either
on the space of Gibbs measures or on the path space) rather we want to find an approximation
that is constructed with low computational cost and with controlled and computable error estimates. 

In general terms we define the coarse-graining operator $\COP:\SIGMA\to\SIGMAC$,
where the coarse configuration space $\SIGMAC$ is defined on the coarse lattice
$\LATTC$, and with the new state space $\SPINSPC$, i.e., $\SIGMAC = {(\SPINSP^c)}^{\LATTC}$. 
The coarse configuration 
$\eta = \COP\sigma \in \SIGMAC$ is defined on a smaller lattice with $M$ lattice
sites and with the coarse state space $\SPINSPC$ for the new lattice spins $\eta(k)$.
The parameter $q$ defines the coarse-graining ratio. 
The operator $\COP$ induces an operator $\COP_*$ on the space of probability measures
$$
\COP_*:\Pp(\SIGMA)\to\Pp(\SIGMAC)\COMMA\;\;\;\; \mu(\sigma) \mapsto 
       \mu^c(\eta) := \mu\{\sigma\in\SIGMA\SEP \COP\sigma=\eta\}\PERIOD
$$ 

\smallskip\noindent
{\it Ising-type spins.}
To be more specific we analyse
the following case of Ising spin-flip dynamics $\SIGMA=\{0,1\}^{\LATT}$. Each coarse
lattice site $k\in\LATTC$ represents a cube $\CUBE_k$ that contains $q$ sites 
of the microscopic lattice $\LATT$. The projection operator defines the block spin
at the coarse site $k$ 
\begin{equation}\label{blockspin}
(\COP\sigma)(k) := \sum_{x\in \CUBE_k} \sigma(x)\PERIOD
\end{equation}
If the dimension $d$ of the lattice is greater than one we understand $k$ and $x$
as multi-indices $k=(k_1,\dots, k_d)$ and we index the corresponding lattice sites 
in the natural order. Choosing the projection operator in this way defines the coarse
state space as $\SPINSPC = \{0,1,\dots q\}$. 
Given the Markov process $(\PROCMIC,\LOPER)$ with the generator $\LOPER$ 
we obtain a coarse-grained process 
$\PROC{\COP\PROCMICRO}$ which is {\it not}, in general, a Markov process. From the computational
point of view this may cause significant difficulties should sampling of such a process
be implemented on the computer. Therefore we derive an {\it approximating} Markov process
$(\PROCMAC,\bar\LOPER^c)$ which can be easily implemented once its generator is given explicitly.

For the model Ising system the projected generator of the coarse-grained process
$\PROCMAC$  can be evaluated explicitly by rearranging the summations on the lattice
$\LATT$. Given the microscopic state $\sigma$ and corresponding coarse state 
$\eta=\COP\sigma$
\begin{eqnarray}\label{cggenerator}
\LOPER\psi(\COP\sigma) = &&\sum_{k\in\LATTC}\left[\sum_{x\in\CUBE_k} c(x,\sigma) (1 - \sigma(x))\right]
                      \left[ \psi(\eta + \delta_k) - \psi(\eta)\right] + \nonumber \\
                      &&\sum_{k\in\LATTC}\left[\sum_{x\in\CUBE_k} c(x,\sigma)\sigma(x)\right]
                      \left[ \psi(\eta - \delta_k) - \psi(\eta)\right]\PERIOD
\end{eqnarray}
The configuration $\delta_k$ defined on the coarse state space is equal to zero at all sites
except the site $k\in\LATTC$ where it is equal $1$, i.e., $\delta_k(j) = 1$ for $j=k$ and
$=0$ otherwise. We see from the formula \VIZ{cggenerator} that the exact generator for the coarse
process can be written in the form
\begin{equation}\label{cgexactgenerator}
 \LOPER^c\psi(\eta) =  \sum_{k\in\LATTC} c_a(k) 
                      \left[ \psi(\eta + \delta_k) - \psi(\eta)\right] +
                      \sum_{k\in\LATTC} c_d(k)
                      \left[ \psi(\eta - \delta_k) - \psi(\eta)\right]\COMMA
\end{equation}
where the new rates 
\begin{equation} \label{newrates}
c_a(k) = \sum_{x\in \CUBE_k} c(x,\sigma) ( 1- \sigma(x))\COMMA\SPACE\SPACE
c_d(k) = \sum_{x\in \CUBE_k} c(x,\sigma) \sigma(x) \COMMA
\end{equation}
correspond to the adsorption and desorption processes. In this form the rates depend on
the microscopic configuration $\sigma$ and not on the coarse random variable $\COP{\sigma}$. 
Therefore, it is reasonable to propose an approximating Markov process, 
which for the case of desorption/adsorption
is a {\it birth-death} process $\PROCMAC$ defined on the state space $\SPINSPC=\{0,1,\dots q\}$.
%
%
This process is defined by the generator $\bar \LOPER^c$ of the form \VIZ{cgexactgenerator} 
where the rates $c_a$ and $c_d$ are replaced by approximate rates
\begin{equation}\label{coarserates}
  \bar c_a(k,\eta) = d_0 ( q- \eta(k))\COMMA\SPACE\SPACE
  \bar c_d(k,\eta) = d_0 \eta(k) \EXP{-\beta\bar U(k,\eta)}\PERIOD
\end{equation}
For details we refer to \cite{kmv1}.
The new rates have a simple interpretation in terms of fluctuations on each cell: $\bar c_a(k,\eta)$
describes the rate with which the coarse variable $\eta(k)$ is increased by one (i.e., adsorption
of a single particle in the coarse cell $\CUBE_k$) and $\bar c_d(k,\eta)$ defines the rate with which
it is decreased by one (desorption in $\CUBE_k$). The new interaction potential $\bar U(\eta)$
represents the approximation of the original interaction $U(\sigma)$.
\begin{definition}
We define the approximation $\BARU(k,\eta)$ of the potential $U(x,\sigma)$, \VIZ{potentialU},
at the coarse level 
\begin{equation}\label{coarsepotentialU}
\BARU(k,\eta) = \sum_{\substack{l\in\LATTC\\ l\neq k}} \bar J(k,l)\eta(l) + 
                  \bar J(0,0) (\eta(k)-1) - \bar h(k)  \PERIOD
\end{equation}
The coarse-grained interaction potential $\bar J$ is computed as the average of the pair-wise
interactions between microscopic spins between the coarse cells $\CUBE_k$ and $\CUBE_l$
\begin{eqnarray}
&& \bar J(k,l) = \frac{1}{q^2} \sum_{x\in\CUBE_k}\sum_{y\in\CUBE_l} J(x-y)\COMMA\;\;\;
     \mbox{for all $k,l\in\LATTC$, such that $k\neq l$, and}\label{coarseJ1} \\
&& \bar J(k,k) \equiv J(0,0) = \frac{1}{q(q-1)} \sum_{x\in\CUBE_k}
     \sum_{\substack{y\in\CUBE_k\\y\neq x}} J(x-y)\PERIOD
\end{eqnarray}
\end{definition}
The error estimate for the projection follows directly from the assumptions on the
regularity of $J$ (or $V$) \VIZ{AssumpV1}--\VIZ{AssumpV2} and the Taylor expansion of the
potential $J$. We state it as a separate lemma.

\begin{lemma}\label{ErrorJ}
Assume that $J$ satisfies \VIZ{AssumpV1}--\VIZ{AssumpV2} then the coarse-grained
interaction potential $\bar J$ at the coarse-graining level $q$ approximates the potential $J$
with the error
\begin{eqnarray}
&& |J(x-y) - \bar J(k,l)| \leq \frac{1}{L} c_d 
          \sup_{\substack{x'\in \CUBE_k\\ y'\in\CUBE_l}}||\nabla V(x'-y')||\leq 
          \BIGO\left(\frac{q}{L^2}\right)
          \label{errJbar1} \\
&& |J(x-y) - \bar J(0,0)| \leq \frac{1}{L} c_d 
             \sup_{\substack{x',y'\in \CUBE_k\\ y'\neq x'}}||\nabla V(x'-y')||
               \leq \BIGO\left(\frac{q}{L^2}\right) \COMMA \label{errJbar2}
\end{eqnarray}
where $c_d = \max_{k\in\LATTC}\{\mathrm{diam}\,(\CUBE_k)\}$.
\end{lemma}

\noindent{\sc Proof:} Using the properties of the potential $V$, we expand
$V$ into the Taylor series, 
$$
V(z) = V(z') + (z-z') . \nabla V(z') + \BIGO(||z-z'||^2)\PERIOD
$$
Using the definition of $J$,  \VIZ{defJV1} and setting $z=x-y$ and $z'=x'-y'$, where
$x,x'\in\CUBE_k$ and $y,y'\in\CUBE_l$, we have
\begin{eqnarray*}
J(x-y) &=& \frac{1}{q^2} \sum_{x'\in\CUBE_k}\sum_{y'\in\CUBE_l} J(x'-y') +\\
       && + \frac{1}{L q^2}\sum_{x'\in\CUBE_k}\sum_{y'\in\CUBE_l} ((x-y) - (x'-y')).\nabla V(x'-y') \\
       && + \frac{1}{L q^2}\sum_{x'\in\CUBE_k}\sum_{y'\in\CUBE_l}\BIGO\left(||(x-y)-(x'-y')||^2\right)\COMMA
\end{eqnarray*}
and using the estimate 
$||(x-y) - (x'-y')||\leq ||x-x'|| + ||y - y'|| \leq \max\{\mathrm{diam}\,(\CUBE_k)\}$
we obtain \VIZ{errJbar1} in the case $k\neq l$ and similarly for $k=l$.

\smallskip

From Lemma~\ref{ErrorJ} we derive the error bound for the approximation
of the coarse-grained potential $\BARU$. Note that in the definition of $U$
the principle contribution to the summation involves interactions within the 
interaction range $L$ and thus we have the following estimate. 
\begin{corollary}\label{ErrorU}
The microscopic potential $U(x,\sigma)$ is approximated by $\BARU(k,\eta)$,
 with the error
\begin{equation}\label{errestim}
\Delta_{q,N}(\BARU,U)\equiv |\BARU(k,\COP\sigma) - U(x,\sigma)| = \BIGO\left(\frac{q}{L}\right)\COMMA
  \;\;\mbox{for all $x\in\CUBE_k$.}
\end{equation} 
\end{corollary}

Note that this approximation represents the direct projection of the interaction kernel $J$
on the coarse space and the contribution from fine scales are neglected. This procedure
differs from the renormalisation group approach where fluctuations from the fine scales
contribute to the transformed Hamiltonian. However, in the case of finite-range
interaction kernels $J$, treated here, the above projection yields approximation of the
order $O(q/L)^2$ as we discuss in the next section. The coarse interaction Hamiltonian
is then given explicitly in terms of $\bar J$ and $\bar h$ as
\begin{equation}\label{coarseHamiltonian}
\bar H(\eta) = -\frac{1}{2}\sum_{l\in\LATTC}\sum_{k\neq l} \bar J(k,l)\eta(k)\eta(l) -
                \frac{1}{2} \bar J(0,0) \sum_{l\in\LATTC} \eta(l)(\eta(l)-1) +
                \sum_{l\in\LATTC} \bar h(l)\eta(l)\PERIOD
\end{equation}

A direct calculation shows that the invariant measure of the Markov process $\PROCMAC$
generated by $\bar \LOPER^c$ is again a canonical Gibbs measure
$$
\mu^c_{M,q,\beta} (d\eta) = \frac{1}{Z_{M,q,\beta}} 
             \EXP{-\beta \BARH(\eta)} \PRODM_{M,q}(d\eta)\COMMA
$$
where the product measure $\PRODM_{M,q}(d\eta)$ is the coarse-grained prior 
distribution. Note that the prior distribution is altered by coarse-graining 
procedure and different projection operators $\COP$ may yield prior distributions that are computationally
intractable.

\noindent
For example, the coarse-grained prior arising from the uniform microscopic prior
($\rho(0) =$ $\rho(1)=1/2$) is the binomial distribution corresponding to
$q$ independent sites:
$$
\PRODM_{M,q}(d\eta) = \PROD_{k\in\LATTC} \rho^c_q(d\eta(k))\COMMA\SPACE\SPACE
\rho^c_q(\eta(k) = p) = \frac{q!}{p! (q-p)!} \left(\frac{1}{2}\right)^{q}\PERIOD
$$

The condition of detailed balance for $\PROCMAC$
with respect to the measure $\mu_{M, q, \beta}$ is
\begin{eqnarray*}
  \bar c_a(k, \eta)\mu_{M, q, \beta}(\eta) &=& \bar c_d(k, \eta+\delta_k)
    \mu_{M, q, \beta}(\eta+\delta_k)\COMMA\\
  \bar c_d(k, \eta)\mu_{M, q, \beta}(\eta) &=& \bar c_a(k, \eta-\delta_k)
    \mu_{M, q, \beta}(\eta-\delta_k)\PERIOD
\end{eqnarray*}
We only verify the first relation, while the second identity is checked in analogous way.
Using that
$
\BARH(\eta+\delta_k)-\BARH(\eta)=-\BARU(k)
$ and the definitions of the rates \VIZ{coarserates},
we have (assuming without loss of generality, $d_0=1$):
\begin{eqnarray*}
&&\bar c_a(k, \eta)\mu_{M, q, \beta}(\eta)- \bar c_d(k, \eta+\delta_k)
            \mu_{M, q, \beta}(\eta+\delta_k) = \\
&&    (q-\eta(k))\EXP{-\beta\BARH(\eta)} \PRODM_{M,q}(\eta)-(\eta(k)+1) 
            \EXP{-\beta\left(\BARH(\eta+\delta_k)+\BARU(k)\right)}  
            \PRODM_{M,q}(\eta+\delta_k)= \\
&&    \EXP{-\beta\BARH(\eta)}\left\{(q-\eta(k))\PRODM_{M,q}(\eta)-
       (\eta(k)+1)\PRODM_{M,q}(\eta+\delta_k) \right\}= \\
&&    \PROD_{l=1, l \ne k}^m \eta(l) 
    \left\{(q-\eta(k)) \eta(k) - (\eta(k)+1) (\eta(k)+1) \right\}\PERIOD
\end{eqnarray*}
Since $(q-p)\rho_q(p)=(p +1)\rho_q(p +1)$,
for all integers  $0 \le p \le q$, 
the last curly bracket  is
equal to zero, hence  the detailed balance holds.
This calculation shows that due to the specific form of the self-interaction
term $\eta(l)(\eta(l)-1)$ the detailed balance condition is satisfied for
the coarse Hamiltonian \VIZ{coarseHamiltonian} and hence the fluctuations from
microscopic dynamics are properly included into the coarse-grained process.
The coarse-graining  procedure described here 
satisfies basic criteria imposed on an approximating process:

\begin{description}
\item[{\rm (i)}]  {error control on a finite-time interval $[0,T]$}. In particular,
    the derived coarse-grained stochastic process $\PROCMAC$ 
    approximates  a pre-specified observable on a finite-time interval $[0,T]$,
    e.g., \VIZ{blockspin}.  In particular, time-dependent error estimates such as 
    \VIZ{timedependent_errorest}
    can rigorously demonstrate that the process $\PROCMAC$ keeps track of  fluctuations 
    from the microscopic level. Consequently expected values of certain path dependent (global) 
    quantities can be properly estimated.
     We characterise approximation properties of $\PROC{\COP\PROCMICRO}$ 
    by $\PROCMAC$ using a suitable probability metric on the path space. 
\item[{\rm(ii)}] approximation of the invariant (equilibrium) measure.
     The invariant measure $\mu^c_{M,q,\beta}(d\eta)$ for the process 
     $\PROCMAC$ defined on $\SIGMAC$ is close, in a suitable probability metric, 
     to the projection of the microscopic measure $\COP_*(\mu_{N,\beta} (d\sigma))$.
     In particular the
     error estimates in \VIZ{apriori_errorest} below demonstrate that the coarse-grained process
     can preserve the 
     ergodicity properties of the microscopic process within a prescribed tolerance. 
     We also note that the coarse-graining
     modifies the microscopic prior $\PRODM_{N}(d\sigma)$ in \VIZ{gibbs}, 
     yielding the coarse prior $\PRODM_{M,q}(d\eta)$.
\end{description}
If the approximating process follows the basic principles (i) and (ii) we observe
as a result of the error estimates presented here and in \cite{kt}, that both
the transient, as well as the long time dynamics are expected to be
captured accurately
by the coarse-graining. Although this is not a complete proof of a
controlled error for  infinite time, it  constitutes a first rigorous
step in this direction. The approximation properties are
also supported by the numerics presented here and in the
references.

%
%
\section{Probability metrics and information theory tools}

Since we propose  the coarse-grained process $\PROCMAC$ to be only an approximation of
$\PROC{\COP\sigma_t}$ which can be computed in a fast and simple way it is necessary
to define in what sense we evaluate the approximation properties. We
view the approximation in coarse-graining procedure as information loss. Such approach
is naturally connected to the actual computational implementation in the Monte Carlo
algorithm. In this section we give a brief introduction to basic tools of information
theory required in the error analysis. We define the basic notions on a probability
space with the countable state space $\STATESP$ but analogous properties and definitions
hold for the relative entropy of  measures on general probability spaces
(see \cite{DupuisEllis}).
Although the exposition in this section is general we keep the notation consistent with
the previous section. However, the reader may assume that the state space $\STATESP$ does
not necessarily refer to the space of spin configurations.

We consider two probability measures $\pi_1(\sigma)$ and $\pi_2(\sigma)$
on the countable state space $\STATESP$, and we define the relative entropy
\begin{equation}\label{relent}
\RELENT{\pi_1}{\pi_2} = \sum_{\sigma\in\STATESP} \pi_1(\sigma) 
                        \log\frac{\pi_1(\sigma)}{\pi_2(\sigma)}\PERIOD
\end{equation}
Using Jensen's inequality it is not difficult to show that
\begin{eqnarray*}
  \RELENT{\pi_1}{\pi_2} &\geq& 0\;\;\;\mbox{and,} \\
  \RELENT{\pi_1}{\pi_2} & = & 0\;\;\;
           \mbox{ if and only if $\pi_1(\sigma) = \pi_2(\sigma)$ for all $\sigma\in \STATESP$.}
\end{eqnarray*}
Although the above properties of the relative entropy $\RELENT{\pi_1}{\pi_2}$ suggest that
this quantity is a distance between the measures $\pi_1$ and $\pi_2$, it does not define
a true metric since it is not symmetric, i.e.,  $\RELENT{\pi_1}{\pi_2}\neq \RELENT{\pi_2}{\pi_1}$
for all measures $\pi_1$, $\pi_2$. Nevertheless, there is an important inequality that
allows us to use the relative entropy as a tool for estimating distance between two measures
and hence use it for evaluating errors in the coarse-graining procedures. 
Using the relative entropy we can bound the total variation of the measures $\pi_1$ and
$\pi_2$:
\begin{equation}
  \RELENT{\pi_1}{\pi_2} \geq \frac{1}{2} \left(\sum_{\sigma\in \STATESP} 
               |\pi_1(\sigma) - \pi_2(\sigma)| \right)^2 \equiv 
               \frac{1}{2} ||\pi_1 - \pi_2||^2_{\mathrm{TV}}\COMMA
\end{equation}
and hence for any observable $\phi=\phi(\sigma)$ we have the bound
\begin{equation}\label{observ}
  |\EXPECWRT{\pi_1}{\phi(\sigma)} - \EXPECWRT{\pi_2}{\phi(\sigma)}|\leq \sup_\sigma |\phi(\sigma)|
                                             \sqrt{2\RELENT{\pi_1}{\pi_2}}\PERIOD
\end{equation}

The following {\it variational} characterisation of the relative entropy
is useful in the error estimation. Given a bounded function (observable)
$\phi\in L^\infty(\STATESP)$ defined on the state space $\STATESP$ we have the natural
dual pairing with the measures on $\STATESP$
$$
  \SCPROD{\pi}{\phi} = \sum_{\sigma\in\STATESP} \pi(\sigma) \phi(\sigma) \equiv
  \EXPECWRT{\pi}{\phi}\PERIOD
$$
The relative entropy \VIZ{relent} has the variational representation (see 
\cite[pp. 338-339]{KL})
\begin{equation}
\RELENT{\pi_1}{\pi_2} = \sup_{\phi\in L^\infty({\STATESP})}
\left\{ \SCPROD{\pi_1}{\phi} - \log \SCPROD{\pi_2}{\EXP{\phi}}\right\}\PERIOD
\end{equation}
The variational representation is used in the next section to obtain
lower bounds on the relative entropy error of coarse-grained processes.

\smallskip

It is worth mentioning the relation between coarse graining, information
theory and application of the relative entropy in the context of coarse graining.
The information point of view also clearly explains the meaning of the relative
entropy as a tool that estimates the loss of information. 
In information theory one is interested in encoding the random variable $\sigma$
with values in the state space $\STATESP$, and distributed according to the 
probability measure $\pi=\pi(\sigma)$, $\sigma\in\STATESP$. The information
should be encoded using symbols from a $D$-nary alphabet, for example only $0$
and $1$ in the case of the binary alphabet. Suppose that $C_D(\sigma)$ is a code/string
corresponding to the value $\sigma\in\STATESP$. We denote $\ell_D(\sigma)$ the
length of the code needed for the state $\sigma$. Since the information is
carried in the random variable $\sigma$ we have to ask what is the {\it expected length}
of the code required to capture the states of $\sigma$ provided we know the distribution
of $\sigma$. The expected length is given by
\begin{equation}
\EXPECWRT{\pi}{\ell_D(S)} = \sum_{\sigma\in\STATESP} \pi(\sigma)\ell_D(\sigma)\PERIOD
\end{equation}
It can be shown (see \cite{ct}) that the optimal (minimal) expected length
is attained by choosing
\begin{equation}\label{explength}
\bar\ell_D(\sigma) = \log_D \frac{1}{\pi(\sigma)}\PERIOD
\end{equation}
Obviously, to set the optimal length for encoding the states of the random
variable $\sigma$ one needs to know the measure $\pi$. If we assume a wrong distribution
$\omega=\omega(\sigma)$ to define the length of the code 
we obtain the expected length which would not be optimal.
The relative entropy $\RELENT{\pi}{\omega}$ describes the increase of the
length \VIZ{explength} due to using the wrong distribution for
the random variable $\sigma$. In this sense $\RELENT{\pi}{\omega}$ is interpreted as
the increase in descriptive complexity due to ``wrong information''.

This information point of view is applicable to the analysis of coarse-graining
procedures: the spin configurations $\sigma$ are sampled by the Markov chain Monte 
Carlo algorithms and hence samples of a random variable $\sigma$ with large-dimensional
state space are generated. On the coarse level we sample an approximate process
$\PROCMAC$ instead of the exact projection $\PROC{\COP\sigma_t}$ and thus
assuming a wrong measure/distribution for the random variable $\sigma$. Using the relative
entropy for evaluating the approximation properties we estimate the loss
of information arising from using samples of $\PROCMAC$ instead of the exact coarse-grained
process. 

%
%
\section{Error analysis and a priori estimates for coarse-grained processes}
As described in the previous section we construct a new process which only
approximates the projected process $\PROC{\COP{\PROCMICRO}}$.
The approximation properties of such construction are quantified in this
section.

We do not attempt to capture the effect of
fine scales exactly and incorporate them into the coarse model through the renormalisation
group transformation. Instead we construct an approximate process $\PROCMAC$, with
the invariant measure $\mu^c_{M,q,\beta}$.
The first question which needs to be addressed is comparison and an error estimate
for the exactly coarse-grained equilibrium measure, i.e., $\COP_*\mu_{N,\beta}$,
and its approximation $\mu^c_{M,q,\beta}$. We recall that $\COP_*$ is the projection
operator induced by the fine-to-coarse projection of spin variables.

\subsection{Information theory estimates} 
The principal idea proposed in \cite{kv} is to control {\it the specific loss of information}
quantified by the relative entropy $\RELENT{\mu^c_{M,q,\beta}}{\COP_*\mu_{N,\beta}}$ between
the coarse-grain equilibrium measure $\mu^c_{M,q,\beta}$ and the projected equilibrium
measure $\COP_*\mu_{N,\beta}$ of the microscopic process.

\begin{proposition}[\cite{kv}, {\it A priori estimate}]
\begin{eqnarray}\label{apriori_errorest}
  &&\frac{1}{N} \RELENT{\mu^c_{M,q,\beta}}{\COP_*\mu_{N,\beta}} := \\
  &&\frac{1}{N} \sum_{\eta\in\SIGMAC} 
    \log\left(\frac{\mu^c_{M,q,\beta}(\eta)}
                   {\mu_{N,\beta}(\{\sigma\in\SIGMA^{\LATT}\SEP
                         \COP\sigma = \eta\})}\right)
                    \mu^c_{M,q,\beta}(\eta) = O\left(\frac{q}{L}\right)\nonumber\PERIOD
\end{eqnarray}
\end{proposition}

This a priori estimate quantifies the dependence of the information distance, the 
specific relative
entropy $\RELENT{\mu}{\nu}$, in terms of the coarse-graining ratio $q$ and the interaction 
range $L$.

The procedure described in the previous section defines a hierarchy of coarse-grained
algorithms parametrised by $q$. The fully resolved simulations correspond to the microscopic
model $q=1$ while the mean-field approximation is obtained in the case where $q\geq L$, i.e.,
when we coarse-grained beyond the interaction range of the potential. Each level of this
hierarchy introduces an error since some fine-scale fluctuations are neglected. 

For the comparison of the processes $\PROC{\COP\PROCMICRO}$ and $\PROCMAC$ we need
to carry out a similar a priori analysis on the coarse path space $\PATHSPC$, i.e.,
on the space of all right-continuous paths $\eta_t: [0,\infty) \to \SIGMAC$. Above
we have presented estimates for the exact coarse graining $\COP_*\mu_{N,\beta}$
of the invariant measure $\mu_{N,\beta}$ and its approximation $\mu^c_{M,q,\beta}$
computed in terms of the coarse Hamiltonian. In a similar way we treat the measures
on the path space: we denote $\PATHMEAS_{\sigma_0,[0,T]}$ the measure on $\PATHSP$
for the process on the interval $[0,T]$, $\PROCT{\sigma_t}$
with the initial distribution $\sigma_0$. Similarly $\PATHMEAS^c_{\eta_0,[0,T]}$
denotes the measure on the coarse path space $\PATHSPC$. With a slight abuse of notation
we also use $\COP_*\PATHMEAS$ to denote the projection of the measure $\PATHMEAS$ 
on the coarse path space, i.e., the exact coarsening of the measure $\PATHMEAS$.
The fully rigorous analysis on the path space is more involved and we refer
to \cite{kt}. For the sake of completeness we only state the main a priori estimate.
\begin{proposition}[\cite{kt}]\label{ktestimate}
Suppose the process $\PROCT{\eta_t}$, given by the coarse generator
$\bar\LOPER^c$, defines the coarse approximation of the microscopic process
$\PROCT{\sigma_t}$ then for any $q<L$ and $N$, $Mq=N$, the information
loss as $q/L\to 0$ is
\begin{equation}
\label{timedependent_errorest}
  \frac{1}{N} \RELENT{\COP_*\PATHMEAS_{\COP_*\sigma_0,[0,T]}}{\PATHMEAS^c_{\eta_0,[0,T]}}
         = T\, O\left(\frac{q}{L}\right)
\end{equation}
\end{proposition}

\noindent{\sc Remark:} The detailed proof of this information estimate (see \cite{kt})
reveals that no control of fluctuations of the process $\PROCMIC$ is necessary for
the estimate. Consequently the estimate is very robust and, as long as $q/L$ is small,
the approximation by the coarse-graining scheme yields a small error independent
of the potential $V$ or the initial distribution $\sigma_0$. Although the estimate is
for finite times $[0,T]$ only, and grows with $T$, in many
cases the system nucleates a new phase at the initial stage of its evolution and thus
the estimate ensures good approximation of the nucleation phase.

\smallskip

It is worth noticing that the relative entropy estimate clearly demonstrate
limitations of the coarse-graining method since it gives the error of order one
for short-range interactions (the nearest neighbour interaction corresponds to $L=1$). On the
other hand the analysis using the relative entropy (information) distance identifies
the small parameter in the asymptotic expansion of the blocking error, namely the
ratio $q/L$.   

The next estimate provides a lower bound for the loss of information
in terms of coarser observables:
\begin{proposition}[{\it Lower bound}]
Suppose the process $(\PROCT{\eta_t},\bar\LOPER^c)$, defined by the coarse-graining
operator $\COP$ with coarse-graining parameters $M q=N$, is the coarse approximation 
of the microscopic process
$\PROCT{\sigma_t}$. Let $\COP^{M',q'}$ be another coarse-graining operator, such that
$M'\leq M$, $M' q' = M q =N$. Then the following estimate for the invariant microscopic
measure $\mu_{N,\beta}$ and the coarse approximation
$\mu^c_{M,q,\beta}$ holds
\begin{equation}
  \RELENT{\mu^c_{M,q,\beta}}{\COP_*\mu_{N, \beta}} \geq 
  \RELENT{\COPP{M',q'}_*\mu^c_{M,q,\beta}}{\COPP{M',q'}_*\mu_{N, \beta}}\PERIOD
\label{apriori}
\end{equation}
Moreover, on any finite-time interval $[0,T]$ 
\begin{equation}
\label{timedependent_errorestII}
  \RELENT{\COP_*\PATHMEAS_{\COP\sigma_0,[0,T]}}{\PATHMEAS^c_{\eta_0,[0,T]}}\geq
  \RELENT{\COPP{M',q'}_*\PATHMEAS_{\COP\sigma_0,[0,T]}}%
         {\COPP{M',q'}_*\PATHMEAS^c_{\eta_0,[0,T]}}\PERIOD
\end{equation}

\end{proposition}

\smallskip
\noindent
{\sc Proof:}
We first recall  the variational formulation for the relative entropy 
\begin{equation}
\RELENT{\mu}{\nu}=\sup_f\left\{ \int f\,d\mu- \log\int e^f \,d\nu\right\}\COMMA
\end{equation}
where the supremum is over all bounded functions in the space where the measures are defined.
This inequality now readily implies the result since 
\begin{equation}
   \RELENT{\mu}{\nu}\geq
    \sup_{f\circ\COP}\left\{ \int f\circ\COP\,d\mu - 
                          \log\int e^{fo\COP} \,d\nu\right\}
    =\RELENT{\COP_*\mu}{\COP_*\nu}
\end{equation}
where $\COP$ is the projection operator (super-scripts omitted) in the statement of the proposition.

\smallskip

\noindent{\sc Remark:} This estimate provides a lower bound for the loss of information
in terms of coarser observables, hence the condition $M' \le M$ where 
$M'q'=Mq=N$. For instance if $M'=1, q'=N$ the measures 
$\COPP{M',q'}_*\mu^c_{M,q,\beta}$ and $\COPP{M',q'}_*\mu_{N,\beta}$
are the PDFs of the total coverage with respect to the coarse-grained
(essentially mean field with a noise) and the microscopic Gibbs 
states respectively. We characterise
such an estimate as {\it a priori} since the 
bound depends on the exact microscopic process, in
analogy to bounds for 
approximations to PDEs which depend on the Sobolev norm of the exact 
solution, \cite{estep}. 
At first glance it may appear that such an estimate is hard to implement since 
it depends on the exact microscopic MC. However, for relatively small systems where 
microscopic MC can be carried out, the bound (\ref{apriori}) can provide a lower bound on
the loss of information,
as well as a sense on how sharp are the upper bounds
given by a posteriori estimates. More specifically
when $M'$ is small , i.e., $M'=1, 2, 3\dots$
etc., the PDFs can be calculated as a histogram
by MC and subsequently the
relative entropy in the lower bound is straightforward
to compute.

%
%
%
\subsection{Microscopic reconstruction and weak convergence estimates}

In many practical MC simulations the main goal is to estimate averages (expected
values) of specific observables. Therefore it is natural to analyze the weak
approximation properties of the coarse-graining procedure. 
The weak error is defined as the quantity
$e_w\equiv|\EXPECWRT{S}{\psi(\COP\PROCMICRO)} - \EXPECWRT{S}{\psi(\PROCMACRO)}|$,
where the expectation $\EXPECWRT{S}{\cdot}$ is defined for the path
conditioned on the initial configuration $\eta_0=\COP\sigma_0 = S$.
Alternatively we can compare the
microscopic process $\PROCMIC$ with its synthetic process $\PROCGAM$ which is
reconstructed from the coarse process $\PROCMAC$. The weak error is then defined
as $e_w\equiv|\EXPECWRT{S}{\phi(\PROCMICRO)} - \EXPECWRT{S}{\phi(\PROCGAMMA)}|$,
where the expectation $\EXPECWRT{S}{\cdot}$ is now defined for the path
conditioned on the initial configuration $\sigma_0 = S$. Here and in what follows
$\phi$ denotes a test function (observable) on the fine level while $\psi$ is used
for a test function on the coarse level.
Theorem~\ref{weakerrorI} and Corollary~\ref{weakerrorII} quantify the
rate of convergence for the weak error on both levels as $q/L\to 0$. 
We refer to \cite{KS} for error estimates in the weak topology
between microscopic MC algorithms and therein derived  SDE  approximations.

Before we formulate the proposition and proceed with the proof it is worth clarifying
the difficulty of comparing the projected process $\PROC{\COP\PROCMICRO}$ with the 
approximating process $\PROCMAC$. The projection $\COP\PROCMICRO$ of the microscopic process
on the coarse grid does not necessarily define a Markov process. On the other hand
the approximating process $\PROCMAC$ is constructed as a Markov process 
$(\PROCMAC,\bar\LOPER^c)$ with the generator $\bar\LOPER^c$ defined by \VIZ{coarserates}.
To circumvent the technical difficulty the authors in \cite{kt} suggested
to construct an auxiliary process $\PROCGAM$ as an intermediate step in the estimation
of the relative entropy between the processes $\PROCMIC$ and $\PROCMAC$. We adopt the
same strategy in order to make comparison between observables which depend
on Markovian processes $\PROCMIC$ and $\PROCGAM$. The process $\PROCGAM$ can be directly
reconstructed from the coarse-grained process $\PROCMAC$.
Thus we are lead to the definition of the {\it synthetic microscopic (Markov) process} 
$\PROCGAM$ associated with the  process $\PROCMIC$.
\begin{definition}[Synthetic microscopic process]\label{syntheticproc}
The auxiliary process $\PROCGAM$ is defined on the microscopic configuration space $\SIGMA$
by the generator $\LOPER^\gamma: L^\infty(\Ss_N) \to \R$ 
\begin{equation}\label{auxgenerator}
(\LOPER^\gamma\phi)(\sigma) = \sum_{x\in \LATT} c_\gamma(x,\sigma) (\phi(\sigma^x) - \phi(\sigma))
\COMMA
\end{equation}
where the rate function $c_\gamma(x,\sigma)$ is defined in terms of the coarse-grained
interaction potential 
$$
c_\gamma(x,\sigma) = d_0(1-\sigma(x)) + d_0\sigma(x)\EXP{-\beta \BARU(k(x),\COP\sigma)}
\PERIOD
$$
The coarse-grained interaction potential $\BARU(k,\eta)$ has been defined in 
\VIZ{coarsepotentialU}. The piece-wise constant interpolation is used
to extend the function $\BARU(.,.)$ from the coarse
lattice to the fine lattice. We denote $k(x)$ to be the cell index of the cell to which
the site $x$ belongs, i.e., $x\in\CUBE_{k(x)}$. 
\end{definition}

The properties of $\PROCGAM$ were studied in
\cite{kt} and it was proved that:
\begin{description}
\item[{\rm (i)}] the coarse-grained projection $\PROC{\COP\PROCGAMMA}$ of
                 the Markov process $(\PROCGAM,\LOPER^\gamma)$ is still a Markov process.
\item[{\rm (ii)}] the processes $\PROC{\COP\PROCGAMMA}$ and $\PROCMAC$ have the same 
                  transition rates. Hence, whenever the processes have the same initial distribution
                  they induce the same probability measure on the coarse-grained path space
                  $\PATHSPC$. If we define $Q^c_{\eta_0}(\eta,t)$ and 
                  $Q_{\gamma_0}(\gamma,t)$
                  to be the probability measures of the Markov processes $\PROCMAC$ and $\PROCGAM$
                  respectively (conditioned on the initial condition $\eta_0 = \COP\gamma_0$), then
                  for all $t>0$ we have the projection
                  $$
                    Q^c_{\eta_0}(\eta,t) = \COP_*Q_{\gamma_0}(\gamma,t) \equiv 
                        \sum_{\{\gamma\SEP\COP\gamma = \eta_t\}} Q_{\gamma_0}(\gamma,t)\COMMA
                  $$
                  provided this relation is satisfied at $t=0$. Hence this property allows us to compare
                  the processes in a path-wise way.
\item[{\rm (iii)}] the microscopic process $\PROCGAM$ can be reconstructed from the approximating
                   coarse process $\PROCMAC$. Such reconstruction is an inverse procedure
                   to the projection from fine to coarse configuration space. In such a way we can
                   compare the original microscopic process with the approximation on the
                   coarse configuration space. A simple choice of a reconstruction operator
                   is to distribute spins $\gamma_t(x)$ for $x\in\CUBE_k$ uniformly so that
                   $\COP\gamma_t|_{\CUBE_k} = \eta_t(k)$.
\end{description}

\noindent{\sc Remark:} It is conceivable that the synthetic process $\PROCGAM$ can be used
not only as a technical tool but as a systematic procedure for reconstructing the
microscopic process $\PROCMIC$ for the purpose of model refinement or adaptivity 
since, as shown in Theorem~\ref{weakerrorI}, the reconstruction is done under rigorous
error estimates. 
In the estimates derived below we deal with a specific class of test functions
$\phi\in L^\infty(\SIGMA)$ which depend
only on the coarse variable $\eta=\COP\sigma$, in other words we impose the assumption
\begin{eqnarray}\label{assumptionA1}
\mbox{(A1)}\;\;\;\;\;\;\;\;
   && \phi(\sigma) = \psi(\COP\sigma)\COMMA\;\;\;\mbox{where $\psi\in L^\infty(\SIGMAC)$, and}
      \label{A1a}\\
   && \sum_{x\in \LATT} |\partial_x \phi(\sigma)| \leq C\COMMA\;\;\;
      \mbox{where $C$ is a constant independent of $N$.}\label{A1b}
\end{eqnarray}

\noindent{\sc Remark:} Observables, such as, for example, the total coverage, used in the numerical 
simulations satisfy this assumption.

\medskip

The principal tool for analysing the weak error is its representation in terms of solutions
to the final value problem on $\SIGMA$
$$
\partial_t v(t,\sigma) + \LOPER v(t,\sigma) = 0\COMMA\;\;\;\;\; v(T,.)= \phi(.) \COMMA
\;\;\; \mbox{for $t< T$} \COMMA
$$
where $\LOPER$ is a generator of the Markov semigroup that defines the lattice dynamics.
Before we state the main estimate of the weak error and its proof we need several preliminary
lemmata that characterize properties of the semigroup generated by the operator $\LOPER$
defined by \VIZ{generator}. 
The specific calculations are better presented by introducing an alternative notation
for the generator $\LOPER$. We define an operator of discrete differentiation
for functions $f\in L^\infty(\SIGMA)$
\begin{equation}\label{discderiv}
\partial_x f(\sigma) \equiv f(\sigma^x) - f(\sigma)\COMMA\;\;\;\mbox{for all $x\in\LATT$,}
\end{equation}
and we introduce two vectors indexed by the lattice sites $x\in\LATT$
$$
\nabla_\sigma f(\sigma) \equiv \left(\partial_x f(\sigma)\right)_{x\in\LATT}\COMMA\;\;\;
\CVEC(\sigma) \equiv \left( c(x,\sigma)\right)_{x\in\LATT}\PERIOD
$$
The scalar product is defined in the natural way as $\CVEC(\sigma)\cdot\nabla_\sigma f(\sigma)
\equiv \sum_{x\in\LATT} c(x,\sigma)\partial_x f(\sigma)$.
Using this notation we write
\begin{equation}\label{generatorII}
\LOPER f(\sigma) = \CVEC(\sigma) \cdot \nabla_\sigma f(\sigma)\COMMA\;\;\;
   \mbox{for all $\sigma\in\SIGMA$.}
\end{equation}
The space of functions defined on the configuration space $\SIGMA$
is equipped with the strong $L^\infty$ topology given by the norm 
$\SUPN{f} \equiv \sup_\sigma\{f(\sigma)\}$.

To prove the estimate in Theorem~\ref{weakerrorI} we need an estimate for the 
difference operator $\nabla_\sigma$ stated here as a separate lemma.
%
%
\begin{lemma}\label{lemma1}
Let $v(t,\sigma)$ be the solution of
\begin{equation}\label{evoleq}
\partial_t v + \LOPER v = 0\COMMA\;\;v(T,\sigma) = \phi(\sigma)\COMMA\;\;
          \mbox{for $t<T$,}
\end{equation}\label{derivestimate}
on a given interval $t\leq T$, then
\begin{equation}\label{bernsteinestimate}
\sum_{x\in\LATT} \SUPN{\partial_x v(t,.)} \leq C_T 
                 \sum_{x\in\LATT} \SUPN{\partial_x \phi} \PERIOD
\end{equation}
Moreover, the constant $C_T$ depends exponentially on the final time $T$.
\end{lemma} 

%
%
\noindent{\sc Proof:} 
Using the notation introduced above
and the definition of $\LOPER$ we recast the evolution equation \VIZ{evoleq} into a familiar
form of a transport equation on the configuration space
\begin{equation}\label{eq1}
\partial_t v + \CVEC(\sigma)\cdot\nabla_\sigma v = 0\COMMA\;\;\;\sigma\in\SIGMA\COMMA\; t>0\PERIOD
\end{equation}
Subtracting \VIZ{eq1} for $v(t,\sigma^x)$ and  $v(t,\sigma)$ we have
$$
\partial_t (v(t,\sigma^x) - v(t,\sigma)) + \CVEC(\sigma)\cdot\left(\nabla_\sigma v(t,\sigma^x) - 
\nabla_\sigma v(t,\sigma)\right) + \left(\CVEC(\sigma^x) - \CVEC(\sigma)\right)\cdot 
\nabla_\sigma v(t,\sigma^x) = 0\COMMA
$$
which we write as
\begin{equation}\label{eq3}
\partial_t\left(\partial_x v(t,\sigma)\right) + \CVEC(\sigma)\cdot
          \nabla_\sigma\left(\partial_x v(t,\sigma)\right) + 
          \partial_x\CVEC(\sigma)\cdot\nabla_\sigma v(t,\sigma^x) = 0\PERIOD
\end{equation}
Next we derive $L^\infty$-bounds for the discrete derivatives $\partial_x\CVEC(\sigma)$
using the explicit definition of the rates $c(x,\sigma)$ in \VIZ{arrhenius_rate}. For each component,
indexed by $z\in\LATT$, of the vector $\CVEC(\sigma)$ we have
$$
\partial_x c(z,\sigma) = c(z,\sigma^x) - c(z,\sigma) = 
  (1-\sigma^x(z)) + \sigma^x(z) \EXP{-U(z,\sigma^x)} - 
  (1-\sigma(z)) + \sigma(z) \EXP{-U(z,\sigma)} \PERIOD
$$
For the spin-flip dynamics, i.e., $\sigma^x(y)=1-\sigma(y)$ if $x=y$ and $\sigma^x(y)=\sigma(y)$
otherwise, a straightforward calculation gives
$\partial_x U(z,\sigma)\equiv U(z,\sigma^x) - U(z,\sigma) = J(z-x) (1-2\sigma(x))$ 
if $z\neq x$ and it is
equal zero otherwise. Thus the discrete derivate $\partial_x\CVEC(\sigma)$ is 
$$
\partial_x c(z,\sigma) = \left\{
           \begin{array}{ll}
                (2\sigma(x)-1)(1-\EXP{-U(x,\sigma)})\COMMA & \mbox{for $z=x$,}\\
                \sigma(z) \EXP{-U(z,\sigma)}\left(1-\EXP{J(x-z)(1-2\sigma(x))}\right) & \mbox{if $z\neq x$.}
           \end{array}\right.
$$
Recalling the definition \VIZ{defJV2} of the interaction potential $J$ we have that $J(z-x)\sim 1/L$ for
$|z-x|\leq L$ and $J=0$ otherwise.
Hence we derived $L^\infty$-bounds for the discrete derivative of the rates
\begin{equation}\label{eq4}
\partial_x c(z,\sigma) \sim \left\{
           \begin{array}{ll}
              \BIGO(1)\COMMA & \mbox{for $z=x$,} \\
              \BIGO(1/L)\COMMA & \mbox{for $|z-x|<L$,} \\
              0 \COMMA & \mbox{otherwise.}
           \end{array}\right.
\end{equation}
Going back to the equation \VIZ{eq3} we have for all $x\in\LATT$
\begin{equation}\label{eq5}
\partial_t\left(\partial_x v(t,\sigma)\right) + \LOPER \partial_x v(t,\sigma) + 
\sum_{z\in\LATT}\partial_x c(z,\sigma)\partial_z v(t,\sigma^x) = 0\PERIOD
\end{equation}
The estimates in \VIZ{eq4} imply
\begin{equation}\label{eq6}
\partial_t \partial_x v(t,\sigma) + \LOPER \partial_x v(t,\sigma) + 
  \BIGO(1)\partial_x v(t,\sigma^x) + 
  \BIGO\left(\frac{1}{L}\right)\sum_{\substack{z\in\LATT\\|z-x|\leq L}}
             \partial_z v(t,\sigma^x) = 0\COMMA
\end{equation}
and we have for all $\sigma\in\SIGMA$ the solution formula 
$$
\partial_x v(t,\sigma) = \EXP{t\LOPER} [\partial_x v(0,\sigma)] + 
               \int_t^T \EXP{(s-t)\LOPER}[\BIGO(1)\partial_x v(s,\sigma^x) + 
  \BIGO(1/L)\sum_{|z-x|\leq L}
             \partial_z v(s,\sigma^x)]\,ds\PERIOD
$$
By the contractive property of the semigroup $\EXP{t\LOPER}$
we have  the estimate
\begin{eqnarray*}\label{eq7}
||\partial_x v(t,\cdot)||_\infty &\leq & ||\partial_x v(0,\cdot)||_\infty + 
  \int_t^T  \BIGO(1) ||\partial_x v(s,\cdot)||_\infty\,ds + \\ 
  && \int_t^T \BIGO(1/L)\sum_{|z-x|\leq L}
             ||\partial_z v(s,\cdot)||_\infty \, ds\COMMA
\end{eqnarray*}
for all $x\in\LATT$.
Thus summing over all $x\in\LATT$ we obtain
\begin{eqnarray*}
\sum_{x\in\LATT}&&||\partial_x v(t,\cdot)||_\infty \leq \sum_{x\in\LATT}||\partial_x v(0,\cdot)||_\infty 
  + \\
&& + \int_t^T (\BIGO(1) \sum_{x\in\LATT} ||\partial_x v(s,\cdot)||_\infty + 
  \BIGO(1/L)\sum_{x\in\LATT}\sum_{|z-x|\leq L}
             ||\partial_z v(s,\cdot)||_\infty )\,ds\COMMA
\end{eqnarray*}
where the last double sum in the integrand is bounded by $2L\sum_x ||\partial_x v(s,\cdot)||_\infty$.
Hence by setting $\theta(t) = \sum_x ||\partial_x v(t,\cdot)||_\infty$ we have
$$
\theta(t) \leq \theta(0) + \int_t^T \BIGO(1) \theta(s)\,ds\COMMA
$$
from which, by using Gronwall's inequality, we obtain the bound
$$
\theta(t) \leq \EXP{c (T-t)} \theta(T)\COMMA
$$
which concludes the proof of \VIZ{derivestimate}.
%
%

Next we establish an $L^\infty$-bound for discrete derivatives of solutions generated
by semigroups $\EXP{t\LOPER}$ and $\EXP{t\LOPER^\gamma}$. 
%
%
\begin{lemma}\label{lemma2}
Let $u(t,\sigma)$ be the solution of
$$
\partial_t u + \LOPER u = 0 \COMMA\;\;\;\; u(T,.) = \phi\COMMA\;\;
\mbox{for $t<T$,}
$$
and let $v(t,\sigma)$ solves
$$
\partial_t v + \LOPER^\gamma v = 0 \COMMA\;\;\;\; v(T,.) = \psi\COMMA\;\;
  \mbox{for $t<T$,}
$$
then for any $t\leq T$ the following estimate holds
\begin{equation}
\sum_{x\in\LATT} \SUPN{\partial_x u(t,\cdot) - \partial_x v(t,\cdot)} \leq
          C_1(T) \sum_{x\in\LATT} \SUPN{\partial_x \phi - \partial_x \psi}
         +C_2(T) \left(\frac{q}{L}\right)\PERIOD
\end{equation}
The constants $C_1$ and $C_2$ are independent of $q$ and $L$ but depend
exponentially on the final time $T$.
\end{lemma}

%
%
\noindent{\sc Proof:} We use the same approach and notation 
as in the proof of Lemma~\ref{lemma1}. Subtracting the evolution equations
and defining $w_x(t,\sigma) \equiv \partial_x u(t,\sigma) - \partial_x v(t,\sigma)$,
$\mathbf{w}(t,\sigma) \equiv (w_x(t,\sigma))_{x\in\LATT}$ we have
\begin{eqnarray}
&& \partial_t w_x(t,\sigma) + \LOPER w_x(t,\sigma) + \label{eq4-l2} \\
&& (\CVEC_\gamma(\sigma)- \CVEC(\sigma))\cdot \nabla_\sigma v(t,\sigma^x) + \label{eq4b-l2}\\
&&  +\partial_x \CVEC(\sigma) \cdot \mathbf{w}(t,\sigma^x) + \label{eq4c-l2} \\ 
&& (\partial_x \CVEC(\sigma) - \partial_x\CVEC_\gamma(\sigma))\cdot 
     \nabla_\sigma v(t,\sigma^x) = 0 \label{eq4d-l2} \PERIOD
\end{eqnarray}
From Lemma~\ref{lemma1} we have estimates for the terms involving $\nabla_\sigma v(t,.)$ (notice
that the lemma essentially gives the estimate of $\SUPN{\nabla_\sigma v(t,.)}$). Furthermore, from
the definition of rates $c(x,\sigma)$ and $c_\gamma(x,\sigma)$ direct calculation (similar to
that used in the proof of Lemma~\ref{lemma1}) yields
the estimate
\begin{equation}\label{eq3-l2}
\SUPN{\CVEC-\CVEC_\gamma} = \BIGO\left(\frac{q}{L}\right)\COMMA
\end{equation}
which allows us to control \VIZ{eq4b-l2} and \VIZ{eq4d-l2}. Term \VIZ{eq4c-l2} is treated
in the same way as a similar term in the proof of Lemma~\ref{lemma1}. Hence, for all $x\in\LATT$
we obtain
$$
\partial_t w_x(t,\sigma) + \LOPER w_x(t,\sigma) + \BIGO(1/L) \sum_{|z-x|\leq L} w_x(z,\sigma^x)
  \leq \BIGO(q/L) \SUPN{\partial_x v(t,.)}\PERIOD
$$
Similarly as in the proof of Lemma~\ref{lemma1} we complete the proof by summing over $x\in\LATT$
and applying Gronwall's inequality.
%
%

\smallskip

Since we are comparing the process $\PROCMIC$ with the process $\PROCGAM$, which is defined
only up to the equivalence given by the projection operator $\COP$ we have to establish 
uniqueness of solutions for initial data satisfying the assumption $(A1)$.
\begin{lemma}\label{lemma3}
Let $\phi\in L^\infty(\SIGMA)$, $\psi\in L^\infty(\SIGMAC)$ be test functions
satisfying (A1). 
Assume that $v(t,\gamma)$ is the solution of the final value problem
\begin{equation}\label{fvpgamma}
   \partial_t v + \LOPER^\gamma v = 0\COMMA\;\;\;
   v(T,\gamma) = \phi(\gamma)=\psi(\COP\gamma)\COMMA
\end{equation}
then for all $\gamma$, $\gamma'\in\SIGMA$ such that $\COP\gamma=\COP\gamma'$
\begin{equation}
v(t,\gamma) = v(t,\gamma')\COMMA\;\;\;\;\mbox{for all $t\leq T$.}
\end{equation}
\end{lemma}

\noindent{\sc Proof:} For convenience we write $v(t,\gamma) = v(t,\COP\gamma)$.
Given a configuration $\gamma\in\SIGMA$ we can reconstruct an arbitrary configuration
$\gamma'\in\SIGMA$ such that $\COP\gamma' = \COP\gamma$ by considering a permutation
$\PERM:\LATT\to\LATT$, $\PERM=(\PERM_1,\dots,\PERM_M)$ such that 
$$
\PERM_k: \CUBE_k\to\CUBE_k\COMMA\;\;k=1,\dots,M\PERIOD
$$
The action of $\PERM$ on the configuration space is defined in a natural way
$\gamma' = \GAMPERM$, or equivalently $\gamma'(x) = \gamma(\PERM x)$. Since the
permutation does not change the total spin in the cell we have $\COP\GAMPERM = \COP\gamma$.
Hence we write $v(t,\gamma') = v(t,\GAMPERM)$ and 
$v(T,\GAMPERM) = v(T,\gamma)$ $= \psi(\COP\gamma)$.
It is sufficient to show that
the function $u(t,\gamma)\equiv v(t,\GAMPERM)$ is a solution of \VIZ{fvpgamma}. From the 
    uniqueness of solutions to \VIZ{fvpgamma} we conclude immediately that $u(t,\gamma)=v(t,\gamma)$.
From the definition of the generator $\LOPER^\gamma$ we have
\begin{equation}\label{eq2-l3}
\partial_t v(t,\GAMPERM) + \sum_{k\in\LATTC}\sum_{x\in\CUBE_k} c_\gamma(x,\GAMPERM)
  (v(t,(\GAMPERM)^x) - v(t,\GAMPERM)) = 0\PERIOD
\end{equation}
Recall the definition of the rate $c_\gamma$
$$
c_\gamma(x,\gamma) = d_0 (1-\gamma(x)) + d_0 \gamma(x) \EXP{-\beta\BARU(k(x),\COP\gamma)}\COMMA
$$
and denote $c_\gamma(x,\gamma)$ by $C_\gamma(\gamma(x),k,\COP\gamma)$ to emphasise the dependence
on $\gamma(x)$, $k$, and $\eta=\COP\gamma$ only.
Thus the inner summation in \VIZ{eq2-l3} becomes
\begin{equation}\label{eq3-l3}
\sum_{x\in\CUBE_k} C_\gamma(\GAMPERM,k,\COP\gamma)  (v(t,(\GAMPERM)^x) - v(t,\GAMPERM))\PERIOD 
\end{equation}
On the other hand the definition of spin-flip dynamics leads to
\begin{equation}\label{eq4-5}
(\GAMPERM)^x(z) = \left\{\begin{array}{ll}
                        \gamma(\PERM z) & z\neq x\COMMA \\
                        1 - \gamma(\PERM x) & z = x\COMMA
                         \end{array}
                   \right. \;\mbox{while}\;
\gamma^{(\PERM x)}(\PERM z) = \left\{\begin{array}{ll}
                        \gamma(\PERM z) & z\neq x\COMMA \\
                        1 - \gamma(\PERM x) & z = x\PERIOD
                         \end{array}
                   \right. 
\end{equation}
Hence we obtain
\begin{equation}\label{eq6-l3}
(\GAMPERM)^x(z) = \gamma^{(\PERM x)}(\PERM z) = \left(\gamma^{\PERM x}\circ \PERM\right)(z)\COMMA
\end{equation}
and substituting to the expression \VIZ{eq3-l3} leads to
\begin{eqnarray*}
&&  \sum_{x\in\CUBE_k} C_\gamma(\gamma(\PERM x),k,\COP\gamma)(v(t,(\GAMPERM)^x) - v(t,\GAMPERM)) = \\
&& = \sum_{x\in\CUBE_k} C_\gamma(\gamma(\PERM x),k,\COP\gamma)(v(t,\gamma^{\PERM x}\circ\PERM) - v(t,\GAMPERM)) = \\
&& = \sum_{y\in\CUBE_k}  C_\gamma(\gamma(y),k,\COP\gamma)(v(t,\gamma^{y}\circ\PERM) - v(t,\GAMPERM)) = \\
&& = \sum_{y\in\CUBE_k}  C_\gamma(\gamma(y),k,\COP\gamma)(u(t,\gamma^{y}) - u(t,\gamma))\PERIOD
\end{eqnarray*}
Thus we have shown that 
$$
\partial_t u(t,\gamma) + \sum_{k\in\LATTC}\sum_{x\in\CUBE_k}  c_\gamma(x,\gamma)(u(t,\gamma^{x})-u(t,\gamma))
=0\PERIOD
$$
Recalling the definition of $u(t,\gamma)$ we obtain that $v(t,\GAMPERM)$ also solves \VIZ{fvpgamma}.
The uniqueness of solutions to \VIZ{fvpgamma} implies that $v(t,\GAMPERM)=v(t,\gamma)$ for all $\gamma$
or $v(t,\gamma') = v(t,\gamma)$ for all $\gamma'$ such that $\COP\gamma'=\COP\gamma$.

\medskip

Now we can formulate and prove the weak error estimate that allows us to compare
the microscopic process and its coarse-level approximation. We estimate the 
weak error on the microscopic level by comparing
the microscopic process and its synthetic process. 
%
\begin{theorem}[{\it Weak error}]\label{weakerrorI}
Let $\phi\in L^\infty(\SIGMA)$ be a test function (observable)
on the microscopic space satisfying (A1) and let 
$(\PROCGAM,\LOPER^\gamma)$ be the synthetic Markov process (in the sense of 
Definition~\ref{syntheticproc}) of the
microscopic process $(\PROCMIC, \LOPER)$ with the initial condition $\sigma_0=S$,
then the weak error satisfies, for $0<T<\infty$, 
\begin{equation}
  |\EXPECWRT{S}{\phi(\sigma_T)} - \EXPECWRT{S}{\phi(\gamma_T)}|
  \leq C_T \left(\frac{q}{L}\right)^2\COMMA
\label{weakII}
\end{equation}
where the constant $C_T$ is independent of $q$ and $L$ but depends on $T$.
\end{theorem}

\medskip
\noindent
{\sc Proof:}
The two ingredients of the proof, the Feynman-Kac formula and
the martingale property, follow from the standard properties
of Markov processes (see for example \cite{KL}). If we define, for the
microscopic process $\PROCMIC$ defined by the generator $\LOPER$, the
function
$$
u(t,S) = \EXPEC{\phi(\sigma_T)\SEP\sigma_t = S}\COMMA
$$
then from the Feynman-Kac formula with the zero potential follows
that the function $u(t,S)$ solves the final value problem
\begin{equation}\label{finvalueprobU}
\partial_t u + \LOPER u = 0\COMMA\;\;\;\;\;u(T,.) = \phi\COMMA\;\;
t<T  \PERIOD
\end{equation}

On the other hand the martingale property implies that for any
smooth function $v(t,S)$ and the process $\PROCGAM$ with the generator
$\LOPER^\gamma$ we have
$$
\EXPECWRT{S}{v(T,\gamma_T)} = \EXPECWRT{S}{v(0,\gamma_0)} + 
   \int_0^T  \EXPECWRT{S}{(\partial_s + \LOPER^\gamma) v(s,\gamma_s)} \,ds\PERIOD
$$
The definition of $u(t,S)$ leads to the representation of the error
$|\EXPECWRT{S}{\phi(\sigma_T)} - \EXPECWRT{S}{\phi(\gamma_T)}|$ by
$e_w = |\EXPECWRT{S}{u(0,S)} - \EXPECWRT{S}{u(T,\gamma_T)}|$
and hence 
$$
e_w =  \left|\int_0^T \EXPECWRT{S}{\left(\partial_s + \LOPER^\gamma\right)u(s,\gamma_s)} \,ds\right|
\PERIOD
$$
The function $u(t,S)$ solves the equation $\partial_t u = -\LOPER u$ thus
we obtain
\begin{eqnarray*}
&&\EXPECWRT{S}{\phi(\sigma_T)-\phi(\gamma_T)}=\int_0^T \EXPECWRT{S}{\LOPER^\gamma u(t,\gamma_t) - \LOPER u(t,\gamma_t)} dt=\\
&=&  \int_0^T  \EXPECWRT{S}{\sum_{x\in \LATT} \left(c(x,\gamma_t) - c_\gamma(x ,\gamma_t)\right)
                           \partial_x u(t,\gamma_t)} dt\PERIOD
\end{eqnarray*}
We split the summation $\sum_{x\in\LATT}$ which gives us
\begin{eqnarray*}
&&\EXPECWRT{S}{\phi(\sigma_T)-\phi(\gamma_T)}= \int_0^T \EXPECWRT{S}{\sum_{k \in{\LATTC}} \sum_{x\in\CUBE_k} 
                                           (c(x,\gamma_t) -c_\gamma(x ,\gamma_t))\partial_x u(t,\gamma_t)} dt=\\
&=&\int_0^T  \EXPECWRT{S}{\sum_{k \in{\LATTC}}\sum_{x\in \CUBE_k}\gamma_t(x)(\EXP{ -\beta  U(x,\gamma_t)}- 
                             \EXP{ -\beta  \BARU(k(x),\COP\gamma_t)})
                             (\partial_k v(t,\COP\gamma_t)+ R^{q,L}_{T}(x))} dt\PERIOD
\end{eqnarray*}
Here we need to replace $\partial_x u$ by the $\partial_x v$, where $v$ solves the final value problem
\VIZ{finvalueprobU} with $\LOPER$ replaced by $\LOPER^\gamma$. From Lemma~\ref{lemma2} we know that the error
term $R^{q,L}_{T}(x) = \partial_x u(t,\gamma) - \partial_x v(t,\gamma)$ 
is controlled by $\BIGO(q/L)$ in $\SUPN{\cdot}$. Furthermore, Lemma~\ref{lemma3} guarantees that
with the final condition $\phi$ which satisfies Assumption (A1) the solution depends only on $\COP\gamma$ and
hence we can replace the discrete difference $\partial_x v$ by the difference $\partial_k v(t,\eta) \equiv
v(t,\eta+\delta_k) - v(t,\eta)$, where $\eta=\COP\gamma$. 
Next  we expand the exponentials to obtain
$$
\Gamma(k,\gamma) \equiv \sum_{x\in\CUBE_k} \beta\gamma(x)\EXP{-\beta\BARU(k(x),\COP\gamma)}
                           \left( \Delta(\BARU,U) + \frac{1}{2} \beta^2\Delta^2(\BARU,U)
                           +\BIGO\left(\beta^3\Delta^3(\BARU,U)\right)\right)\COMMA
$$
and we recast the error representation into
\begin{eqnarray}
&&\EXPECWRT{S}{\phi(\sigma_T)-\phi(\gamma_T)}= 
             \int_0^T  \EXPECWRT{S}{\sum_{k \in{\LATTC}}\Gamma(k,\gamma_t)
                           \partial_k v(t,\COP\gamma_t) + 
            \sum_{x\in\LATT}(c(x,\gamma_t)-c_\gamma(x,\gamma_t)) R^{q,L}_{T}(x)} \, dt \nonumber \\
&=&\int_0^T \EXPECWRT{S}{q \sum_{k \in\LATTC}\partial_k v(t,\eta_t)\,
            \EXPEC{\Gamma(k,\gamma)\big| \COP\gamma=\eta_t}} \,dt \label{expect1} + \\ 
&+& \int_0^T \EXPECWRT{S}{\sum_{x\in\LATT}  (c(x,\gamma_t)-c_\gamma(x,\gamma_t)) R^{q,L}_{T}(x)}\,dt 
             \PERIOD \label{expect2} 
\end{eqnarray}
Assumption~(A1) and Lemma~\ref{lemma1} imply that the term $q \sum_{k \in\LATTC}\partial_k v(t,\eta_t)$
is bounded. To estimate the conditional expectation
we use the property of the reconstruction operator for the process $\PROCGAM$, in particular on each cell
$\PROCGAMMA(x)$ is reconstructed from $\PROCMACRO(k)$ by assuming a ``local'' equilibrium and distributing
$\PROCGAMMA(x)$ uniformly in the cell $\CUBE_{k(x)}$. Using this property we can compute the conditional
expectation explicitly and we obtain for  $l \ne k$
$$
\EXPEC{\sum_{x \in \CUBE_k} \gamma(x) \Delta(\BARU,U)\big| 
   \COP\gamma=\eta}
  =\eta_k\eta_l \sum_{\substack{x \in \CUBE_k\\y \in \CUBE_l}} 
   \left(J(x-y)-\bar J_{kl}\right) = 0\PERIOD
$$
Similarly we handle the case $l=k$ and we conclude that, after averaging, the first-order 
term $\Delta(\bar U,U)$ in $\Gamma(k,\gamma)$ vanishes.
We recall (see \VIZ{ErrorU}) that
$$
\Delta(\BARU,U) \equiv \BARU(k(x),\COP\gamma) - U(x,\gamma) = \BIGO\left(\SMALLPAR\right)\COMMA
$$
and hence we can estimate \VIZ{expect1} by $O(q^2/L^2)$. For the term
\VIZ{expect2} we use 
the estimate  $\sum_{x\in\LATT}|R^{q,L}_{T}(x)|\sim \BIGO(q/L)$ from Lemma~\ref{lemma2} 
and the H\"older inequality
$$
\EXPECWRT{S}{\sum_{x\in\LATT} (c(x,\gamma_t)-c_\gamma(x,\gamma_t)) R^{q,L}_{T}(x)}\leq 
  \SUPN{\mathbf{c}-\mathbf{c}_\gamma}\EXPECWRT{S}{\sum_{x\in\LATT}|R^{q,L}_T(x)|}
\PERIOD
$$
The first term on the right-hand side is estimated from \VIZ{eq3-l2} by $C (q/L)$ and hence the
left-hand side behaves as $\BIGO(q^2/L^2)$. Combining the estimates of \VIZ{expect1} and \VIZ{expect2} 
we conclude the proof.  
%
%

\medskip

%
%
Using the estimate for the synthetic process and its reconstruction from the coarse-grained process
$\PROCMAC$ we can compare the projected process 
$\{\COP\sigma_t\}_{t\geq 0}$ and the coarse-grained process
$\PROCMAC$ also on the coarse level. 
The weak error for observables on the coarse space
is also natural in simulations where we usually project finer simulations on the coarse level
and use estimators for the coarse processes.
\begin{corollary}\label{weakerrorII}
Let $\psi\in L^\infty(\SIGMAC)$ be  a test function on the coarse level
such that there exists a test function $\phi\in L^\infty(\SIGMA)$ 
satisfying (A1) with the property
$\psi(\COP\sigma)=\phi(\sigma)$. Given the initial configuration
$\sigma_0$ we define the coarse configuration $\eta_0=\COP\sigma_0$. Assume the 
microscopic process $(\PROCMIC, \LOPER)$ with the initial condition 
$\sigma_0$ and 
the approximating coarse process $(\PROCMAC,\bar\LOPER^c)$ with the
initial condition $\eta_0=\COP\sigma_0$, then the weak
error satisfies, for $0<T<\infty$, 
\begin{equation}
  |\EXPECWRT{S}{\psi(\COP\sigma_T)} - \EXPECWRT{S}{\psi(\eta_T)}|
  \leq C_T \left(\frac{q}{L}\right)^2\COMMA
\label{weak}
\end{equation}
where the constant $C_T$ is independent of $q$ and $L$ but depends on $T$.
\end{corollary}

%
%
\section{Implementation of the coarse-grained Monte Carlo algorithms}
The hierarchy of coarse-grained Monte Carlo processes (CGMC) parametrised by $q$
has been designed in such a way that it is easily implemented in the unified manner.
In fact, the nature of the generator $\bar \LOPER^c$ at the level $q$ allows us to use
the same implementation as for the standard MC at the microscopic level, i.e., $q=1$.

The stochastic system is simulated with the kinetic Monte Carlo (KMC) algorithm. 
Each iteration of  the Monte Carlo simulation produces a variable time step $\Delta t$ within which 
a spin flip occurs at a specific lattice node based on the transition probability,
$$
[c_a(k,\eta) + c_d(k, \eta)] \Delta t + O(\Delta t^2)
$$
where $c_a$ and $c_d$ are as in (\ref{coarserates}). This procedure repeats until the stopping 
criteria (see below) have been met. More specifically, the simulation is implementing the following 
{\it global updating} process-type kinetic Monte Carlo (KMC) algorithm for spin 
flip Arrhenius dynamics:
\begin{description}
\item[{\rm Step 1}] Calculate all transition rates $c_a(k,\eta)$ (adsorption), 
      $c_d(k,\eta)$ (desorption) from (\ref{coarserates}) for all nodes $k$ 
      in the lattice $\LATTC$
\item[{\rm Step 2}] Calculate the total $R_a = \sum_{l\in \LATTC} c_a(l, \eta)$, 
      $R_d = \sum_{l \in \LATTC} c_d(l,\eta)$ 
      adsorption, desorption rates respectively. 
      Similarly obtain the total rate $R_T = R_a + R_d$.
\item[{\rm Step 3}] Obtain two random numbers $\rho_1$ and $\rho_2$.
\item[{\rm Step 4}] Use the first random number to choose between absorption or desorption 
      based on the measure created by the rates $R_a, R_d$ and $R_T$. 
      Assume that the choice is to adsorb(desorb) and  denote by 
      $c \equiv c_a(l,\eta)$, $(c_d(l,\eta))$ and $R = R_a$, $(R_d)$, respectively.
\item[{\rm Step 5}] Find the node at lattice position $l \in\LATTC$ such that,
$$
\sum_{j=0}^l c(j,\eta) \geq \rho_2 R \geq \sum_{j=0}^{l-1} c(j,\eta)
$$
\item[{\rm Step 6}] Update the time, $t = t + \Delta t$ where 
\begin{equation}
\Delta t = 1/R_T\PERIOD
\label{timest}
\end{equation}
\item[{\rm Step 7}] Repeat from Step 1 until equilibrium or dynamics of interest have been captured.
\end{description}

As expected a kinetic Monte Carlo algorithm produces no ``null'' steps and therefore 
every trial is accepted. A similar version of the algorithm can also be implemented 
with a {\it local} updating mechanism which can improve speed substantially at 
the reciprocal expense of allocating further computer memory for dynamic array allocation. 
In the simulation that follow we use a finite size interaction potential and lattice 
size $L, N<\infty$.

We produce simulations and compare observables at microscopic ($q=1$) and coarse grained ($q>1$) levels. 
For consistency purposes we use the same seed for our random number generator in order to compare 
simulations for different coarse grained values of $q$. This allows us to focus on the differences 
attributed only to the coarse graining variable and not on those resulting from different paths due 
to the initial seed. In the case of several realisations we initialise each new microscopic realisation 
with a different seed. Once again, for comparison purposes, we  initialise each subsequent 
coarse grained realisation with the same seeds used in the respective microscopic simulations.
All simulations are compared in the same non-dimensional time units. The corresponding non-dimensional 
time-step is respectively set by the Monte Carlo simulation based on the rule \ref{timest}.

In the simulations which follow we try, when possible, to group together various parameters in the model
so the results are presented with respect of variations in a smaller number of parameters. 
In that respect we point out that for the fixed external field $\bar{h}$ it is possible to 
group together $d_0$ and $\bar{h}$ in (\ref{coarserates}) as follows,
$$
\bar{c}_d (k,\eta) = c_0 \eta (k) e^{-\beta [\sum_l \bar{J}(k,l) \eta(l) + \bar{J}(0,0) (\eta(k) - 1)]} 
$$
where $c_0 = d_0 e^{\bar{h}}$. We provide the values of all pertinent parameters as well as $d_0$ and $c_0$ 
in the relevant figures.

%
%
\section{Numerical simulations}

We use the CGMC described and analyzed in the previous sections for efficient simulations
in the spin systems that undergo phase transitions. 
Within the context of spin-flip dynamics a typical example is
nucleation of spatial regions of a new phase or a transition from one phase (all spins equal to zero)
to another (all spins equal to one). In such simulations the emphasis is on the path-wise properties of
the coarse-grained process so that the switching mechanism is simulated efficiently while approximation
errors are controlled. 
We compare simulations on the microscopic level $q=1$ with
those performed on different levels of coarse-graining hierarchy parametrized by $q$.

\begin{figure}
\centerline{\hbox{\psfig{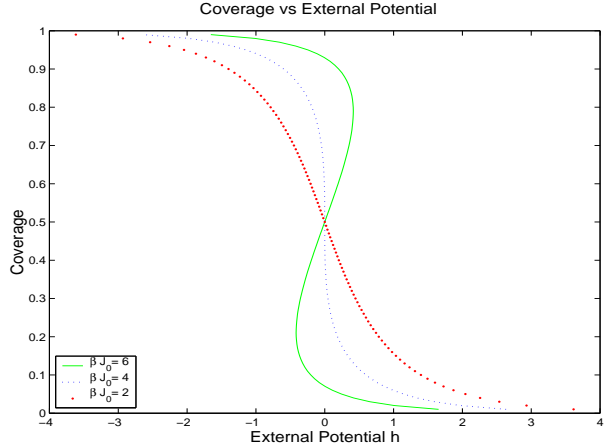}}}
\caption{Equilibrium coverage $c(\sigma)$ and its dependence on the external field. The critical
         point for the $\{0,1\}$ spins satisfies $\beta_c J_0 = 4$. The solid line depicts equilibria
         below the critical temperature. The hysteresis shape of the curve manifests existence
         of two equilibria in the neighbourhood of the zero external field.}
\label{hysteresis}
\end{figure}

The qualitative behaviour of the Ising model with a long-range potential can be understood
from the mean-field approximation of the equilibrium total coverage $c(\sigma)$. Below the critical
temperature the Gibbs measure is not unique (in the thermodynamic limit $N\to\infty$) and two phases 
can coexist. When the energy landscape
is probed by changing the external field $h$ we observe non-uniqueness of the equilibrium coverage
as depicted in Figure~{\ref{hysteresis}}. The fluctuations allow for transitions between the 
equilibrium which leads to nucleation of regions with a different phase. Changing the external field
$h$ makes the original phase unstable and a switching occurs -- the system transforms
into the other equilibrium configuration.

%

The parameters in the simulations have been chosen as follows:
We use a uniform finite range potential for all examples presented. We simulate a finite lattice 
with a total of $N=1000$ microscopic nodes and allow a potential interaction range of 
$2 L + 1$ for $L=100$. 
We choose the constant $d_0 = 1$ so that $c_a = 1$ and $c_d = 1$. 
Hence in this case the critical value $\beta_c$ is given by $\beta_c J_0 =  4$.
If $\beta J_0 > \beta_c J_0 = 4$ the system is in the phase transition regime
and the two phases can coexist. In this region we typically observe a transition
from one phase (e.g., zero (low) coverage) to the other phase (e.g., full coverage).  
For the phase transition examples we fix $\beta J_0 = 6 > \beta_c J_0$. The simulations
become difficult when $\beta \simeq \beta_c$ and there is no external field $h$ applied.
We note that the coarse-graining algorithm will not perform well close to the critical
point $\beta_c$ when $h=0$. 
In the numerical studies we first investigate approximation properties of the CGMC 
algorithms for certain global quantities. 

\smallskip\noindent
{\it Coverage:} We define the coverage
$c_t$ to be the process computed as the spatial mean
$$
  c_t(\sigma_t)      = \frac{1}{N} \sum_{x\in\LATT}\sigma_t(x)\COMMA\SPACE\SPACE
  c^{q}_t(\eta_t)  = \frac{1}{qM} \sum_{l\in\LATTC}\eta_t(k)\PERIOD
$$

%
\begin{figure}
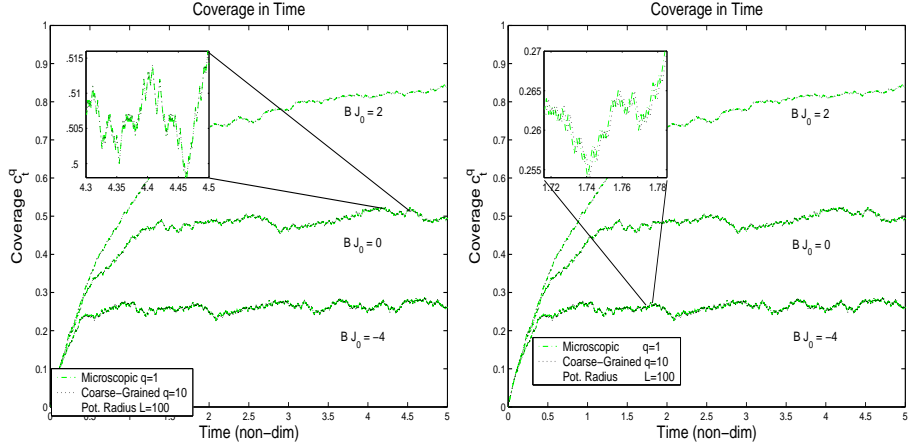

\centerline{\hbox{\psfig{figure=Figures/relaxationcompmicrocg.epsc,height=6cm,width=6cm}}
            \hbox{\psfig{figure=Figures/relaxationcompmicrocg2.epsc,height=6cm,width=6cm}}}
\caption{Relaxation dynamics. Comparison of microscopic ($q=1$) and  coarse grained ($q=10$) 
         simulations. The plot depicts a short time simulation in order to calibrate the code 
         and compare to Figure 4 from \cite{kmv1}.}
\label{Timeseries2}
\end{figure}
%
%
%
\begin{figure}
\centerline{\hbox{\psfig{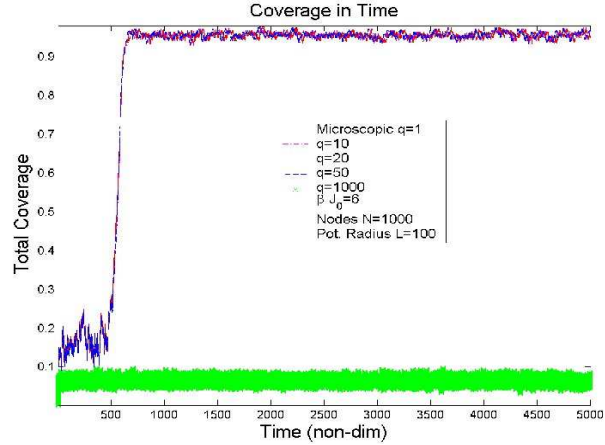}}}
\caption{Time series of the coverage $c^q_t$. 
         Simulations for different coarse-graining ratios are shown in the phase transition regime. 
         The case $q=1000, m=1$ (mean-field approximation) shows significant discrepancy.  
         Parameters used: potential radius length $L=100$, $\beta J_0 = 6$, $d_0 = 1$, $c_0 = .072$}
\label{Timeseries}
\end{figure}

We present time evolution of the coverage at the phase transition regime, $\beta J_0 = 6$. 
Note that the case $q=1000$, $m=1$ which corresponds to the mean-field approximation 
(``over coarse-grained'' interactions)
does not follow the phase transition path of the other simulations. 
On the other hand  the agreement in the results is extremely good for the remaining values of $q$.
Furthermore, these numerical results indicate path-wise (strong) approximation of the microscopic process
by the coarse-grained process. This observation suggests a stronger error control than the
relative entropy estimate provided by Proposition~\ref{ktestimate}.

To quantify the error behaviour we calculate two errors between the exact stochastic
process $c_t$ and its coarse approximation $c^q_t$ at the level of coarse-graining $q$.
We define the weak error  $e_w[c]$ and the strong error $e_s[c]$ respectively:
$$
e_w[c] = \int_0^T \left|\EXPEC{c_t} - \EXPEC{c^q_t}\right| \d t\COMMA\;\;\;\;
e_s[c] = \int_0^T \EXPEC{|\COP c_t - c^q_t|} \d t\PERIOD
$$
The expected values are estimated by empirical means and the integral in time by the
piece-wise constant quadrature.

\begin{figure}
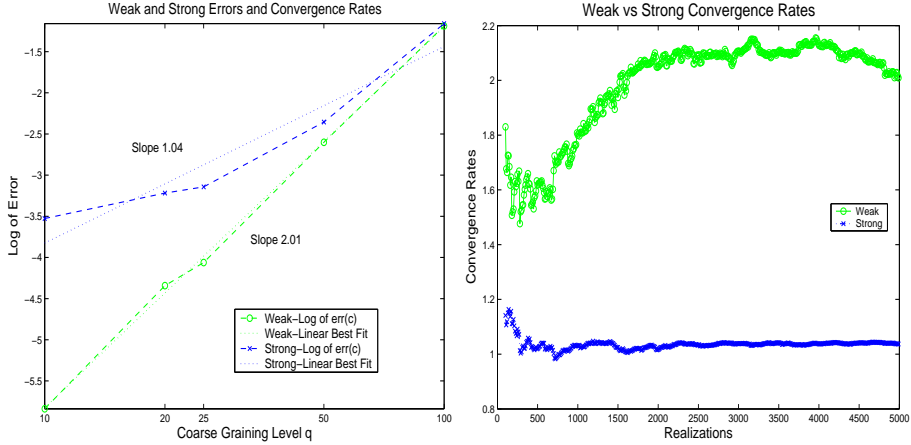

\centerline{\hbox{\psfig{figure=Figures/weakvsstronglogofslope.epsc,height=6cm,width=6cm}}
            \hbox{\psfig{figure=Figures/weakvsstrongslopeestim.epsc,height=6cm,width=6cm}}}
\caption{Estimated weak $e_w[c]$ and strong $e_s[c]$ errors. We compare the exact process
         $c_t$, $q=1$ with coarse approximations $c^q_t$, $q=10,25,50$ and $100$. The simulation
         parameters were fixed at $L=100$, $d_0=1$, $c_0=.07$, $\beta J_0 = 6 > \beta_c J_0$ and 
         the lattice size $N=1000$. The convergence rates depicted are estimated by the linear
         best fit on the logarithmic scale. The statistical error or dependence of the estimates
         on the number of realisations is depicted in the right figure.
        }
\label{figslope}
\end{figure}

The simulations allow us to estimate the convergence rate for both errors. The rates in the case
of fixed parameters $L=100$, $d_0=1.0$, $c_0=0.07$ and $\beta J_0 = 6$ on the lattice of the size
$N=1000$ are depicted in Figure~\ref{figslope}. Note that we need
to eliminate the statistical error, arising from approximation of expected values by empirical means.
However, as seen in Figure~\ref{figslope} the estimator of the rate converges as the number of realisations
tends to infinity.

%
%
%

Since the coarse-grained Hamiltonian neglects higher order corrections arising
from the fluctuations on fine scales one may expect that the approximation
is poor if $q/L$ is not very small. This is certainly true at the critical point
(i.e., $\beta = \beta_c$ and $h=0$) but further from the critical point the 
approximation properties are improved. This is demonstrated in the following table, where
the simulations were performed in the presence of different (large) external
fields. The relative error becomes small even for fairly crude coarse-graining $q=20$
in the case of shorter interaction radii $L$.

\begin{table}
\begin{center}
\label{errctable}
\caption{Relative strong error $e_s[c]$ in the presence
of an external field defined by $c_0$. Comparisons are made for different
values of the interaction radius $L$ and different coarse-graining
levels $q$. Size of the lattice fixed at $N=1000$.}
\begin{tabular}{|l||l|lll|}
\hline
   $c_0$    & $L$       &   $q=5$   &  $q=10$ & $q=20$ \\ 
\hline 
            & 100       &     .0591 &  .0733  &  .1134 \\
.07         & 40        &     .0820 &  .0880  &  .1113 \\
            & 20        &     .1508 &  .2214  &  .1832 \\
\hline
            & 100       & .0186     & .0563   &  .0480 \\
.09         &  40       & .0678     & .0749   &  .1064 \\
            &  20       & .1760     & .1767   &  .1812 \\
\hline
            & 100       & .0010     & .0010   &  .0025  \\
1           & 40        & .0036     & .0040   &  .0054  \\
            & 20        & .0016     & .0043   &  .0065  \\
\hline
\end{tabular}
\end{center}
\end{table}

\smallskip\noindent
{\it Mean time to reach phase transition:}
One quantity of interest that is calculated from the simulations is the mean time 
$\bar\tau_T=\EXPEC{\tau_T}$ until the coverage reaches 
$C^{+}$ in its phase transition regime (see Figure \ref{Timeseries}).
The  random exit time is defined as $\tau_T = \inf\{t>0\SEP c_t\geq C^+\}$. We estimate
the probability distributions $\rho_\tau$ and $\rho_\tau^q$ from the simulations.
We record a phase transition at the time $\bar{\tau}_{T}$ when the coverage exceeds the threshold value $C^+=0.9$.
%
%
\input{ttable.tex}
%

\begin{figure}[h]
\centerline{
\hbox{\psfig{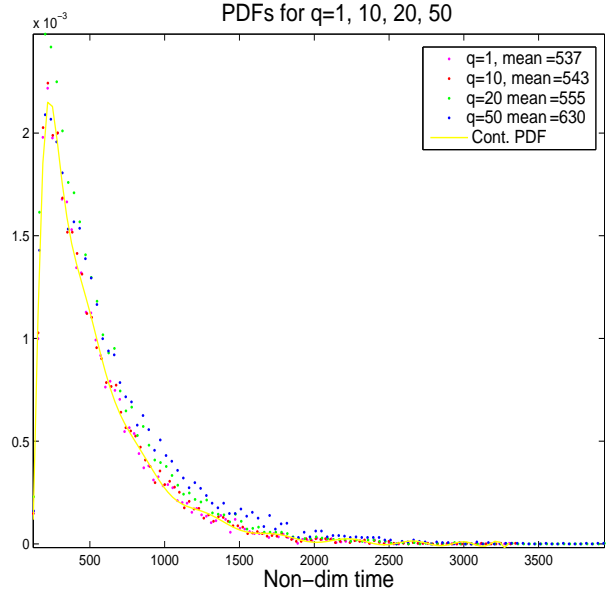}}}
\caption{Probability Density Function (PDFs) comparisons between different coarse-graining 
levels $q$. 
The estimated mean times for each PDF are shown in the figures. 
All PDFs comprised of 10000 samples and the histogram is approximated by 100 bins.}
\label{pdfs}
\end{figure}

In Figure~\ref{pdfs} we plot approximations of the Probability Density Functions (PDFs) 
of $\tau_T$ and compare them for different values of $q$. 
%
The qualitative agreement observed in Figure~\ref{pdfs} is  quantified by using the information
distance for error estimation, i.e., by estimating the relative entropy
\begin{equation}
\RELENT{\rho_1}{\rho_2} = \sum_\lambda \rho_1(\lambda) 
                          \log\left(\frac{\rho_1(\lambda)}{\rho_2(\lambda)}\right)\PERIOD
\label{relativentropy}
\end{equation}

\begin{figure}
\centerline{
            \hbox{\psfig{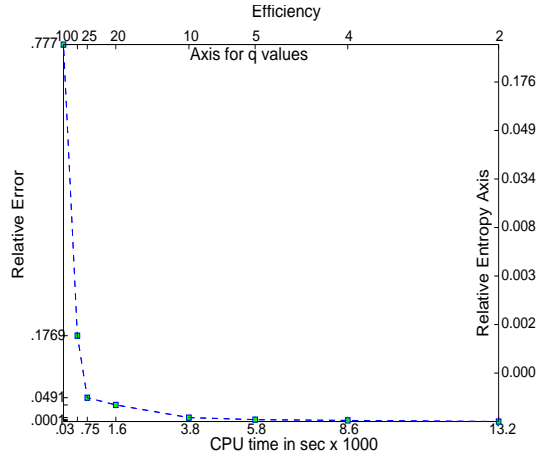}}
}
\caption{The dependence of the relative error and the relative entropy on the coarse-graining level $q$.
         The left scale on the vertical axis depicts the relative error and $q$ on the log-log scale.
         Measurements based on averaging over $10000$ realisations 
         for each $q$.
         }
\label{mixed}
\end{figure}

\smallskip\noindent
{\it Nucleation:} The nucleation of a new phase is a typical phenomenon in the regime
where $\beta > \beta_c$. Essentially, there exist two equilibria (phases).
Random fluctuations will induce transitions from one state to another by overcoming 
energy barriers that separate the equilibria.
We investigate approximation of the path-wise behaviour 
on the configuration space for nucleation of a new phase. 
Two different initial configurations are used. 

\begin{figure}[h]
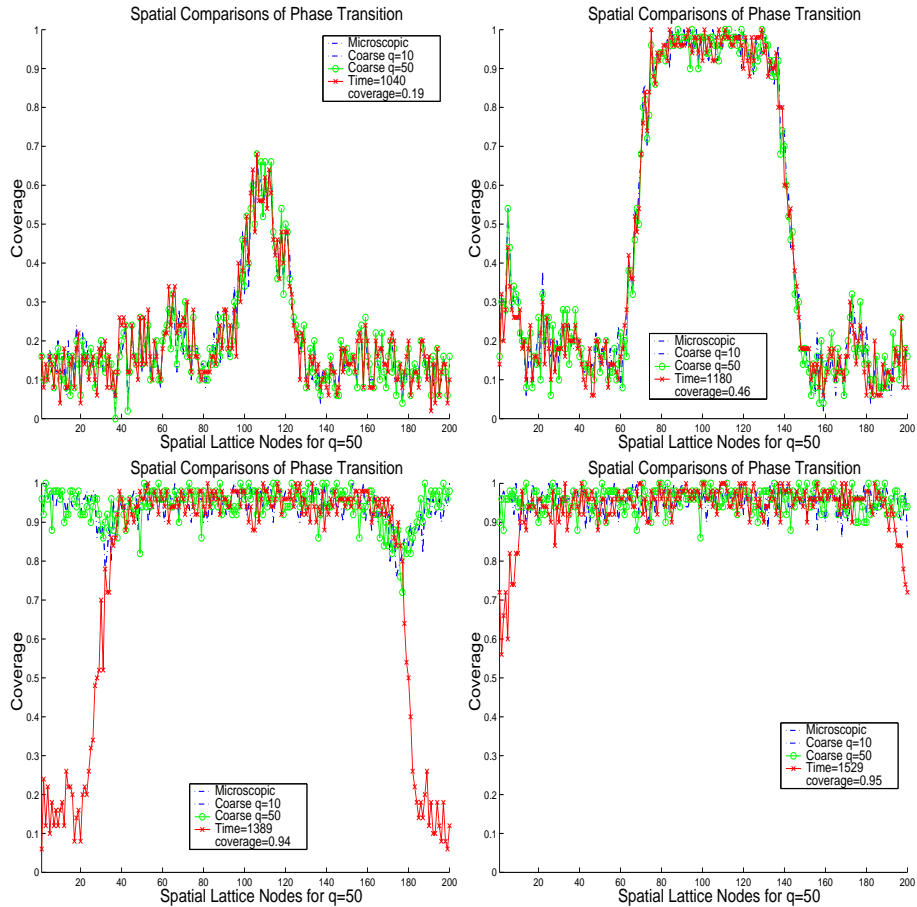

\centerline{
            \hbox{\psfig{figure=NoseedL200/noseed3.epsc,height=6cm,width=6cm}}
            \hbox{\psfig{figure=NoseedL200/noseed5.epsc,height=6cm,width=6cm}}}
\centerline{
            \hbox{\psfig{figure=NoseedL200/noseed8.epsc,height=6cm,width=6cm}}
            \hbox{\psfig{figure=NoseedL200/noseed10.epsc,height=6cm,width=6cm}}}
\caption{Snap-shots of the transition from zero initial spatial distribution. 
         Comparisons between the microscopic $q=1$ and two coarse grained 
         simulations $q=10$ and $q=50$. The interaction radius is set to $L=200$ while 
         total nodes are $N=10000$.}
\label{noseedL200}
\end{figure}

\noindent {\sc Test Case I:} The initial state is at the metastable equilibrium where
the coverage is zero. 
The fluctuations will cause the transition to the full coverage equilibrium which
is stable due to the external applied field.
We present only qualitative comparison in the series
of snap-shots (Figure~\ref{noseedL200}) of the phase transition from 
the uniform (zero) 
initial coverage to the full coverage. We observe a striking path-wise agreement on
the configuration space for relatively large values of $q$ compared to the interaction
radius $L$. However, as the ratio $q/L$ increases the corresponding coarse-grained process
lags behind which is also demonstrated in the expected values of transition times.
Such behaviour suggests that fluctuations at regions with uniform states are well-approximated
by a highly coarse-grained process while finer resolution is necessary for resolving nucleation
of new phases through islands. 

\begin{figure}[h]
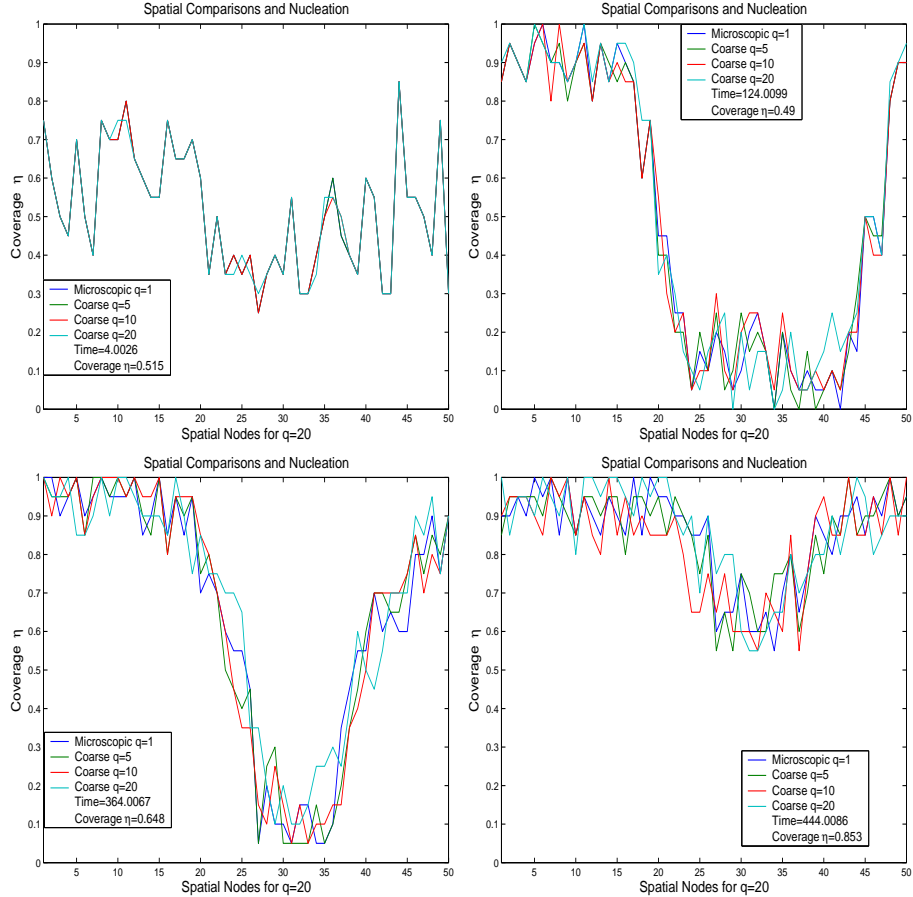

\centerline{
            \hbox{\psfig{figure=SaddleL100/Covp5SpatFig1.epsc,height=6cm,width=6cm}}
            \hbox{\psfig{figure=SaddleL100/Covp5SpatFig3.epsc,height=6cm,width=6cm}}}
\centerline{
            \hbox{\psfig{figure=SaddleL100/Covp5SpatFig5.epsc,height=6cm,width=6cm}}
            \hbox{\psfig{figure=SaddleL100/Covp5SpatFig6.epsc,height=6cm,width=6cm}}}
\caption{Snap-shots of the transition from the initial state with the mean coverage at $0.5$. 
         Comparisons between the microscopic $q=1$ and coarse grained 
         simulations $q=5, 10$ and $q=20$. The interaction radius is set to $L=100$,
         the external field  $c_0=0.0492$, $d_0=1$ and the total number of lattice 
         sites $N=1000$.}
\label{SaddleL100}
\end{figure}

\noindent {\sc Test Case II:} 
 We have already documented the path-wise agreement of the approximating dynamics 
under both transition and relaxation cases. In this example we examine the nucleation 
phenomenon at the critical parameter regime of phase 
transition $\beta J_0 = 6$. We chose the initial state to be at a saddle point of the energy
surface, i.e., the mean coverage is set to $0.5$. Snapshots of the spatial distribution 
of spins are presented in Figure \ref{SaddleL100}. Under all four dynamics examined 
$q=1, 5, 10$ and $20$ we observe complete spatial path-wise 
agreement. Over time the total coverage may fall towards zero or rise towards one 
in which case it will remain 
there since we are at the phase transition regime where these represent stable equilibria. 
Furthermore such spatial examples of nucleation are shown below in Test Case III under 
the assumption of an ``island-type'' of initial state.


\begin{figure}
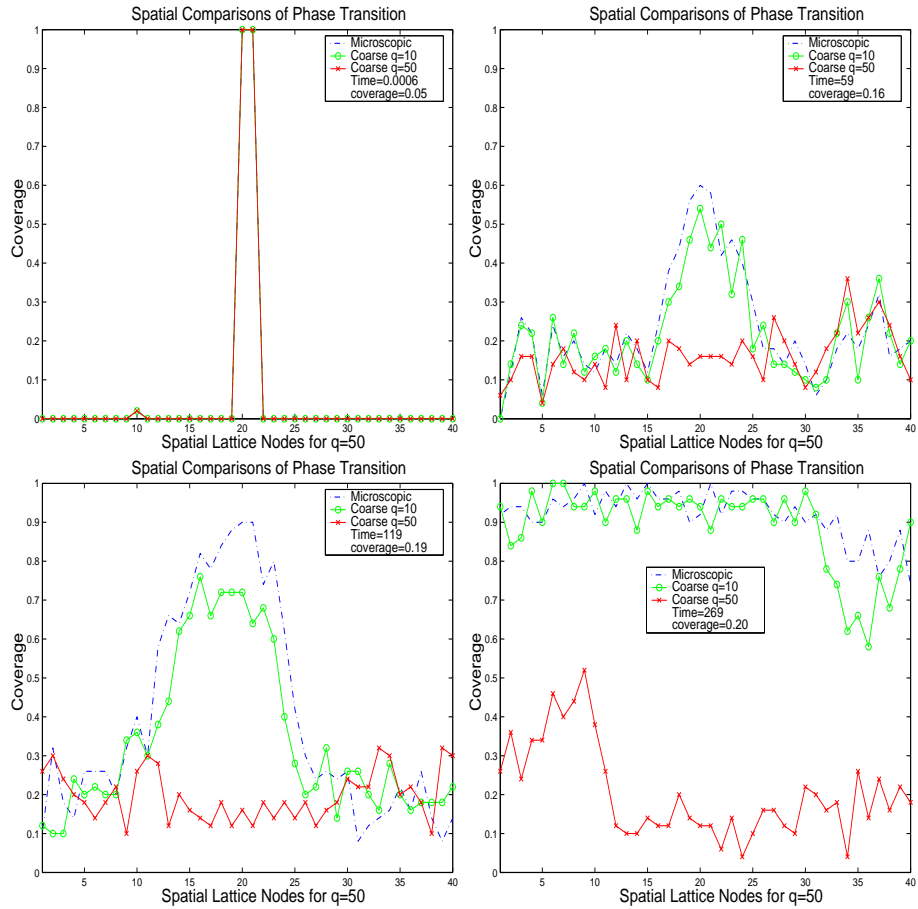

\centerline{\hbox{\psfig{figure=Smallseed/smallseed1.epsc ,height=6cm,width=6cm}}
            \hbox{\psfig{figure=Smallseed/smallseed3.epsc ,height=6cm,width=6cm}}}
\centerline{\hbox{\psfig{figure=Smallseed/smallseed5.epsc ,height=6cm,width=6cm}}
            \hbox{\psfig{figure=Smallseed/smallseed10.epsc,height=6cm,width=6cm}}}
\caption{Snap-shots of the nucleation from a small-size initial seed. 
         Comparisons between the microscopic $q=1$ and two coarse 
         grained simulations $q=10$ and $q=50$. Potential radius is set to $L=100$ and 
         the size of the lattice to $N=2000$.}
\label{smallseed}
\end{figure}

\noindent {\sc Test Case III:} The last set of simulations presents evolution
from the non-uniform initial state, giving a qualitative comparison of nucleation from an island
of a given size (Figure~\ref{smallseed}).
In these simulations we observe spatial propagation of the interface 
in time for different initial size
of the island.

%

\medskip\noindent
{\bf Acknowledgements:}
The research of M.K. is partially supported by NSF-DMS-0413864 and
NSF-ITR-0219211. The research of P.P. was partially supported by 
NSF-DMS-0303565.  
M.K. also acknowledges many valuable conversations with A. Szepessy.
The authors would like  to thank
the Institute for Mathematics and its Applications
where part of this work was carried out during the programme 
``Mathematics of Materials and Macromolecules: Multiple Scales, Disorder, 
and Singularities''


\bibliographystyle{siam}
\bibliography{cgmc}


\end{document}

%% file: ppmacro1.tex
\def\LATT{{\Lambda}_N}
\def\LATTC{{\Lambda^{c}_{M}}}

\def\SPINSP{\Sigma}
\def\SPINSPC{\Sigma^c}
\def\CONFMICRO{\sigma}

\def\PROCMICRO{\sigma_t}
\def\PROCMACRO{\eta_t}
\def\PROCGAMMA{\gamma_t}

\def\PROCMIC{\{\sigma_t\}_{t\geq 0}}
\def\PROCMAC{\{\eta_t\}_{t\geq 0}}
\def\PROCGAM{\{\gamma_t\}_{t\geq 0}}
\def\PROC#1{\{{#1}\}_{t\geq 0}}
\def\PROCT#1{\{{#1}\}_{t\in[0,T]}}
\def\JNOT{V}
\def\PRODM{\mathrm{P}}
\def\STATESP{\mathcal{S}}
\def\SIGMA{{\mathcal{S}_N}}
\def\SIGMAC{{\mathcal{S}^c_{M,q}}}
\def\CUBE{C}
\def\COP{\mathbf{T}}
\def\COPP#1{\mathbf{T}^{#1}}

\def\CVEC{\mathbf{c}}
\def\SUPN#1{{||{#1}||_\infty}}

\def\PERM{\pi}
\def\GAMPERM{\gamma\circ\PERM}

\def\LOPER{\mathcal{L}}
\def\BARH{\bar H}
\def\BARU{\bar U}

\def\RELENT#1#2{\mathcal{R}\left({#1}\SEP{#2}\right)}

\def\PATHSPC{\mathcal{D}(\SIGMAC)}
\def\PATHSP{\mathcal{D}(\SIGMA)}
\def\PATHMEAS{Q}

\def\Pp{{\mathcal{P}}}

\def\Pp{\mathcal{P}}

\def\Ss{\mathcal{S}}

\def\R{\mathbb{R}}

\def\T{\mathbb{T}}
\def\S{\mathbb{S}}
\def\Z{\mathbb{Z}}

\def\d{d}

\def\DT{\Delta t}

\def\EXP#1{e^{#1}}

\def\EXPEC#1{{\mathbb{E}\left[{#1}\right]}}

\def\EXPECWRT#1#2{{\mathbb{E}_{#1}\left[{#2}\right]}}

\def\SMALLPAR{\frac{q}{L}}

\def\SPACE{\;\;\;}          
\def\COMMA{\,,}             
\def\PERIOD{\,.}            
\def\SEP{{\,|\,}}           

\def\VIZ#1{(\ref{#1})}      

\def\BIGO{\mathrm{O}}

\def\SCPROD#1#2{\langle {#1},{#2}\rangle}

\def\PROD{\prod}

\def\Sph{{\S}}


%% file: ttable.tex
\begin{table}
\begin{center}
\label{table21}
\caption{Approximation of $\bar\tau_T, \RELENT{\rho^q_\tau}{\COP_*{\rho_\tau}}$ and relative error.}
\begin{tabular}{|rrrrrr|}\hline
   $L$          &     $q$     &   $\bar\tau_T$   &  $\RELENT{\rho^q_\tau}{\COP_*{\rho_\tau}}$ &\mbox{Rel. Err.}& CPU [s] \\ \hline
   100          &      1      &   532         &          0.0           &0       &  309647  \\
   100          &      2      &   532         &          0.003         &0.01\%  &  132143   \\
   100          &      4      &   530         &          0.001         &0.22\%  &  86449   \\ 
   100          &      5      &   534         &          0.003         &0.38\%  &  58412  \\
   100          &      10     &   536         &          0.004         &0.82\%  &  38344  \\
   100          &      20     &   550         &          0.007         &3.42\%  &  16215  \\
   100          &      25     &   558         &          0.010         &4.91\%  &  7574  \\
   100          &      50     &   626         &          0.009         &17.69\% &  4577  \\
   100          &      100    &   945         &          0.087         &77.73\% &  345  \\
\hline
\end{tabular}
\end{center}
\end{table}

%% file: MAKPPAS-kmc.bbl
\begin{thebibliography}{10}

\bibitem{BRA}
{\sc D.~R. A.~Brandt and D.~J. Amit}, {\em Multi-level approaches to
  discrete-state and stochastic problems}, in Multigrid Methods, II,
  W.~Hackbusch and U.~Trottenberg, eds., Springer-Verlag, 1986, pp.~66--99.

\bibitem{Bai}
{\sc D.~Bai and A.~Brandt}, {\em Multiscale computation of polymer models.}, in
  Multiscale Computational Methods in Chemistry and Physics, J.~B. A.~Brandt
  and K.~Binder, eds., vol.~177 of NATO Science Series: Computer and System
  Sciences, IOS Press, 2001, pp.~250--266.

\bibitem{Talay-Tubaro1}
{\sc P.~Bernard, D.~Talay, and L.~Tubaro}, {\em Rate of convergence of a
  stochastic particle method for the {K}olmogorov equation with variable
  coefficients}, Math. Comp., 63 (1994), pp.~555--587, S11--S17.

\bibitem{bcr}
{\sc G.~Beylkin, R.~Coifman, and V.~Rokhlin}, {\em Fast wavelet transforms and
  numerical algorithms. {I}}, Comm. Pure Appl. Math., 44 (1991), pp.~141--183.

\bibitem{B92}
{\sc A.~Brandt}, {\em Multigrid methods in lattice field computations}, Nucl.\
  Phys.\ B\/, 26 (1992), pp.~137--180.

\bibitem{BI2}
{\sc A.~Brandt and V.~Ilyin}, {\em Multilevel {M}onte {C}arlo methods for
  studying large-scale phenomena in fluids.}, J.\ of Molecular Liquids, 105
  (2003), pp.~253--256.

\bibitem{BrandtRon:01}
{\sc A.~Brandt and D.~Ron}, {\em Renormalization multigrid: Statistically
  optimal renormalization group flow and coarse-to-fine {M}onte {C}arlo
  acceleration}, J. Stat. Phys., 102 (2001), pp.~231--257.

\bibitem{ct}
{\sc T.~M. Cover and J.~A. Thomas}, {\em Elements of Information Theory}, New
  York: Wiley, 1991.

\bibitem{DupuisEllis}
{\sc P.~Dupuis and R.~S. Ellis}, {\em A weak convergence approach to the theory
  of large deviations}, John Wiley \& Sons Inc., New York, 1997.
\newblock A Wiley-Interscience Publication.

\bibitem{ee}
{\sc W.~E and B.~Engquist}, {\em Multiscale modeling and computation}, Notices
  Amer. Math. Soc., 50 (2003), pp.~1062--1070.

\bibitem{estep}
{\sc K.~Eriksson, D.~Estep, P.~Hansbo, and C.~Johnson}, {\em Computational
  differential equations}, Cambridge University Press, 1996.

\bibitem{Sokal}
{\sc J.~Goodman and A.~D. Sokal}, {\em Multigrid {M}onte {C}arlo methods for
  lattice field theories.}, Phys.\ Rev.\ Lett., 56 (1986), pp.~1015--1018.

\bibitem{golden}
{\sc Q.~Hou, N.~Goldenfeld, and A.~McKane}, {\em Renormalization group and
  perfect operators for stochastic differential equations}, Phys. Rev. E, 6303
  (2001), p.~6125.

\bibitem{hou}
{\sc T.~Y. Hou and X.-H. Wu}, {\em A multiscale finite element method for
  {PDE}s with oscillatory coefficients}, in Numerical treatment of multi-scale
  problems (Kiel, 1997), vol.~70 of Notes Numer. Fluid Mech., Vieweg,
  Braunschweig, 1999, pp.~58--69.

\bibitem{steph}
{\sc A.~E. Ismail, G.~C. Rutledge, , and G.~Stephanopoulos}, {\em
  Multiresolution analysis in statistical mechanics. {I}. using wavelets to
  calculate thermodynamic properties.}, J. Chem. Phys., 118 (2003),
  pp.~4414--4424.

\bibitem{steph1}
{\sc A.~E. Ismail, G.~C. Rutledge, and G.~Stephanopoulos}, {\em Multiresolution
  analysis in statistical mechanics. {II}. wavelet transform as a basis for
  {M}onte {C}arlo simulations on lattices}, J. Chem. Phys., 118 (2003),
  p.~4424.

\bibitem{Kadanoff}
{\sc L.~P. Kadanoff}, {\em Statistical physics : statics, dynamics and
  renormalization}, World Scientific, 1999.

\bibitem{KPRT}
{\sc M.~Katsoulakis, P.~Plech\'a\v{c}, L.~Rey-Bellet, and D.~Tsagkarogiannis},
  {\em Higher order corrections to the coarse-grained {M}onte {C}arlo method.}
\newblock in preparation, 2005.

\bibitem{kms}
{\sc M.~A. Katsoulakis, A.~J. Majda, and A.~Sopasakis}, {\em Multiscale
  couplings in prototype hybrid deterministic/stochastic systems: Part i,
  deterministic closures}, Comm. Math. Sci., 2 (2004).

\bibitem{kmv1}
{\sc M.~A. Katsoulakis, A.~J. Majda, and D.~G. Vlachos}, {\em Coarse-grained
  stochastic processes and {M}onte {C}arlo simulations in lattice systems}, J.
  Comput. Phys., 186 (2003), pp.~250--278.

\bibitem{kmv2}
{\sc M.~A. Katsoulakis, A.~J. Majda, and D.~G. Vlachos}, {\em Coarse-grained
  stochastic processes for microscopic lattice systems}, Proc. Natl. Acad. Sci.
  USA, 100 (2003), pp.~782--787 (electronic).

\bibitem{KS}
{\sc M.~A. Katsoulakis and A.~Szepessy}, {\em Stochastic hydrodynamical limits
  of stochastic spin systems}.
\newblock in preparation, 2005.

\bibitem{kt}
{\sc M.~A. Katsoulakis and J.~Trashorras}, {\em Information loss in
  coarse-graining of stochastic particle dynamics}, tech. report, University of
  Massachussets, 2004.
\newblock accepted in J. Stat. Phys.

\bibitem{kv}
{\sc M.~A. Katsoulakis and D.~G. Vlachos}, {\em Hierarchical kinetic {M}onte
  {C}arlo simulations for diffusion of interacting molecules}, J. Chem. Phys,
  112 (2003).

\bibitem{kevrek}
{\sc I.~G. Kevrekidis, C.~W. Gear, and G.~Hummer}, {\em Equation-free: the
  computer aided analysis of complex multiscale systems}, AIChE Journal, 50
  (2004), pp.~1346--1355.

\bibitem{KL}
{\sc C.~Kipnis and C.~Landim}, {\em Scaling Limits of Interacting Particle
  Systems}, Springer-Verlag, 1999.

\bibitem{Kloeden}
{\sc P.~E. Kloeden and E.~Platen}, {\em Numerical Solution of Stochastic
  Differential Equations}, vol.~23 of Applications of Mathematics, Springer
  Verlag, 3rd~ed., 1999.

\bibitem{bi}
{\sc D.~P. Landau and K.~Binder}, {\em A Guide to {M}onte {C}arlo Simulations
  in Statistical Physics}, Cambridge University Press, 2000.

\bibitem{mtv}
{\sc A.~J. Majda, I.~Timofeyev, and E.~V. Eijnden}, {\em A mathematical
  framework for stochastic climate models}, Comm. Pure Appl. Math., 54 (2001),
  p.~891.

\bibitem{m-p}
{\sc F.~Muller-Plathe}, {\em Coarse-graining in polymer simulation: from the
  atomistic to the mesoscale and back}, Chem. Phys. Chem., 3 (2002), p.~754.

\bibitem{schuette1}
{\sc C.~Sch{\"u}tte, A.~Fischer, W.~Huisinga, and P.~Deuflhard}, {\em A direct
  approach to conformational dynamics based on hybrid {M}onte {C}arlo}, J.
  Comput. Phys., 151 (1999), pp.~146--168.
\newblock Computational molecular biophysics.

\bibitem{Szepessy1}
{\sc A.~Szepessy, R.~Tempone, and G.~E. Zouraris}, {\em Adaptive weak
  approximation of stochastic differential equations}, Comm. Pure Appl. Math.,
  54 (2001), pp.~1169--1214.

\bibitem{Talay-Tubaro2}
{\sc D.~Talay and L.~Tubaro}, {\em Expansion of the global error for numerical
  schemes solving stochastic differential equations}, Stochastic Anal. Appl., 8
  (1990), pp.~483--509 (1991).

\bibitem{vlachos90}
{\sc D.~G. Vlachos, L.~D. Schmidt, and R.~Aris}, {\em The effects of phase
  transitions, surface diffusion and defects on surface catalyzed
  reactions:fluctuations and oscillations}, J. Chem. Phys., 93 (1990), p.~8306.

\end{thebibliography}
